\documentclass[10pt]{article}
\usepackage{amsfonts, amsmath}
\usepackage[pdftex]{graphics}
\usepackage[hang]{subfigure}
\usepackage{rotating}

\renewcommand{\theenumi}{\roman{enumi}}

\makeatletter \renewcommand\p@enumi{(\theenumi)} \makeatother

\def\complaint#1{}
\def\withcomplaints{
\newcounter{mycomplaints}
\def\complaint##1{\refstepcounter{mycomplaints}%
\ifhmode%
\unskip%
{\dimen1=\baselineskip \divide\dimen1 by 2 %
\raise\dimen1\llap{\tiny -\themycomplaints-}}\fi%
\marginpar{\tiny [\themycomplaints]: ##1}}%
}

\newtheorem{proposition}{Proposition}

\newtheorem{example}{Example}
\newtheorem{remark}{Remark}

\usepackage{tikz}
\usetikzlibrary{arrows}

\begin{document}

\title {States of self-stress in symmetric frameworks and applications}

\author{Bernd Schulze\\Department of Mathematics and Statistics, Lancaster University\\ Lancaster, LA1 4YF, UK\\ \\
Cameron Millar\\ Skidmore, Owings \& Merrill\\
The Broadgate Tower, 20 Primrose St, London, EC2A 2EW, UK \\\\
Arek Mazurek\\ Mazurek Consulting\\
1012 Frances Pkwy, Park Ridge, IL 60068, USA \\\\
William Baker \\ Skidmore, Owings \& Merrill\\
224 S. Michigan Avenue, Suite 1000, Chicago, IL 60604, USA
}

\maketitle

\begin{abstract}
We use the symmetry-extended Maxwell rule established by Fowler and Guest to detect states of self-stress in symmetric planar frameworks.  The dimension of the space of self-stresses that are detectable in this way may be expressed in terms of the number of joints and bars that are unshifted by various symmetry operations of the framework. Therefore, this method provides an efficient tool to construct symmetric frameworks with many `fully-symmetric' states of self-stress, or with `anti-symmetric' states of self-stress. Maximizing the number of independent self-stresses of  a planar framework, as well as understanding their symmetry properties, has important practical applications, for example in the design and construction of gridshells. We show the usefulness of our method by applying it to some practical examples. 

\end{abstract}

\section{Introduction} 

This paper investigates states of self-stress and mechanisms in symmetric 2D bar-joint frameworks. Such frameworks consist of pin-jointed nodes and axially rigid members. In the field of mathematical rigidity theory, these frameworks are represented as straight line realisations of \emph{graphs} in the plane. Attention is restricted to \emph{planar} frameworks  in which no two bars cross each other, since these are of particular interest in structural engineering applications. However, the methods also extend to non-planar frameworks in a straightforward fashion.

A key tool in this paper is the symmetry-adapted counting rule for bar-joint frameworks developed by Fowler and Guest \cite{FGsymmax}, which extends  the conventional Maxwell count \cite{Cal}.
The derivation of the Fowler-Guest counting rule relies upon \emph{group theory}. Many practitioners are not familiar with the mathematical theory of groups, but since the resulting rule only involves counting bars and nodes with certain symmetry properties, the method is very quick and easy to use. An accompanying paper aims to give a simplified non-technical description of the Fowler-Guest counting rule and its applications discussed here \cite{Millar_sym}. The present paper focuses on unpinned frameworks, but the methods easily extend to pinned frameworks, as discussed in Section \ref{sec:pinned}. 

Motivation for this paper comes from the design of gridshell structures, such as the Great Court Roof of the British Museum, London. Such structures project down onto the $xy$ plane to produce a \textit{form diagram} \cite{Millar2021}. Millar et al \cite{Millar_sym} discuss the role of the states of self-stress in the form diagram within the design of gridshells. It is desirable for gridshells to be quad-dominant, so the examples in this paper focus on quad-dominant frameworks. Quadrilateral glass panels tend to be cheaper than triangular panels as there is less material wastage in their manufacture.  Furthermore, the nodes can be torsion free (the members at a node share a common axis). This is seldom the case for triangulated gridshells which have many high-valent nodes. 

\begin{figure}[htp]\begin{center} 
      \includegraphics[scale=0.23]{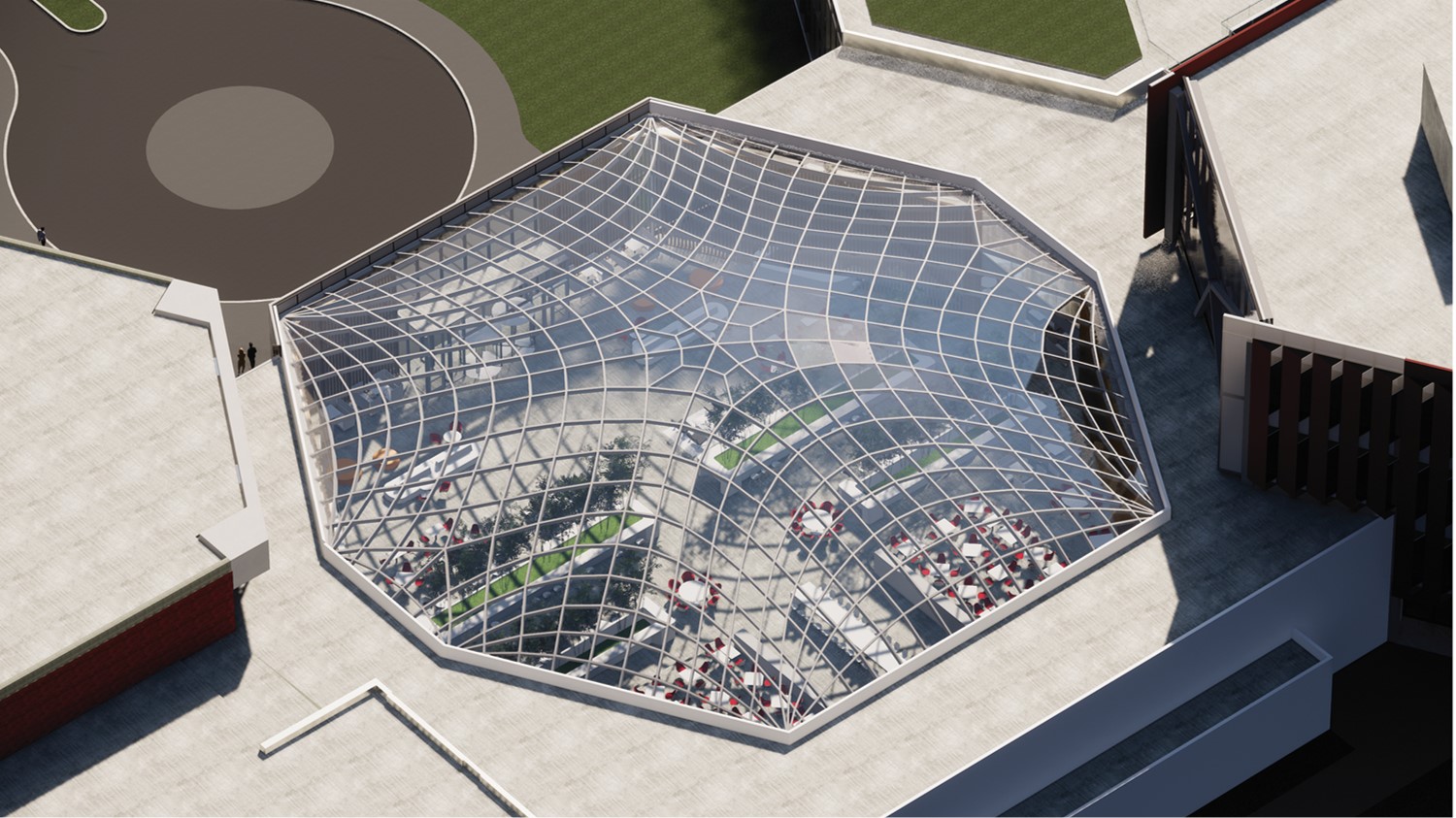}
      \hspace{0.1cm}
       \includegraphics[scale=0.23]{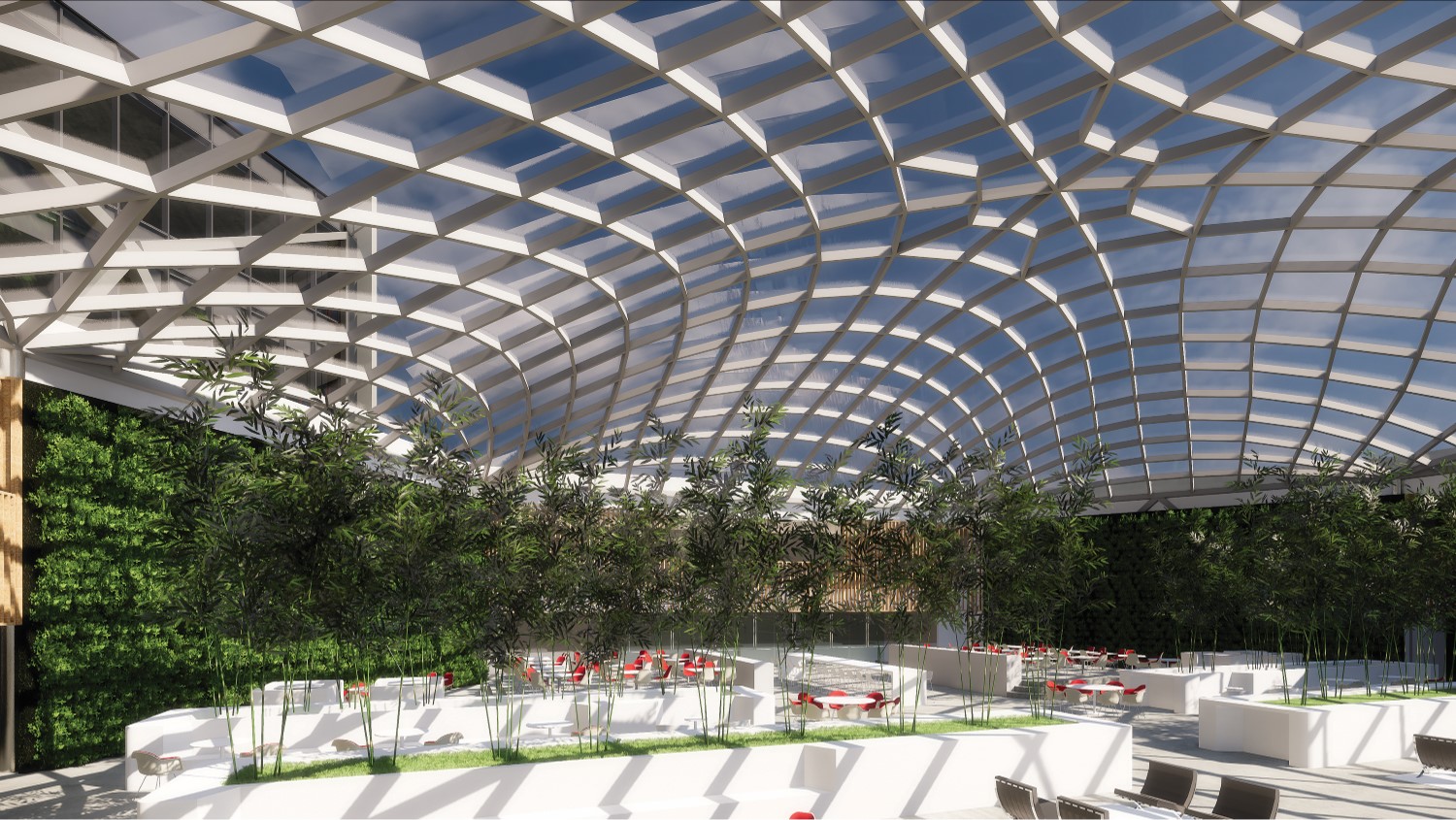}
      \end{center}
      \caption{A  symmetric quad-dominant gridshell structure.}
\label{fig:pic}
\end{figure}

Many gridshell structures possess symmetry (see Figure~\ref{fig:pic} for an example) so it is natural to try and utilise the symmetry-adapted counting rule as an analysis and design tool.  
As the states of self-stress of 2D frameworks are a projectively invariant property \cite{Izmestiev,Nixon}, it is possible to design highly symmetric frameworks with many states of self-stress and then project them to obtain a geometry which fits the construction requirements. Such an example is discussed in Section~\ref{sec:dih}. As noted in Section~\ref{sec:formulas}, using a larger symmetry group can increase the number of states of self-stress detected.

Millar et al \cite{Millar2021} describe how each state of self-stress of the form diagram relates to a funicular gravity loading of the gridshell (the applied loads are taken through axial forces only -- there is no bending moment in the gridshell). Funicularity is a desirable engineering property as it can reduce the volume of material needed to construct the load-bearing gridshell structure. Therefore, one often wants to increase the number of states of self-stress within the form diagram so that the size of the funicular load space is increased accordingly. A fully-symmetric state of self-stress relates to a symmetric vertical loading which is also preferable (self-weight is an important and sometimes dominant load case which is symmetric). Anti-symmetric states of self-stress relate to an anti-symmetric loading of the gridshell. Pattern loading of the gridshell (uneven gravity loads) can often be decomposed into a fully-symmetric and anti-symmetric load, as discussed by McRobie et al \cite{Allan2020}. Therefore, anti-symmetric states of self-stress can be a useful property when designing gridshells. 

\begin{figure}[htp]
\begin{center}
\begin{tikzpicture}[very thick,scale=0.7]
\tikzstyle{every node}=[circle, draw=black, fill=white, inner sep=0pt, minimum width=4pt];
   
       \path (0,0) node (p1)  {} ;
    \path (-1,-1) node (p2)  {} ;
    \path (-2,-1) node (p3)  {} ;
     \path (1,-1) node (p4)  {} ;
      \path (2,-1) node (p5)  {} ;
     \path (0,-1.5) node (p6)  {} ;
     \path (0,-2.5) node (p7)  {} ;
     
               \draw (p1)  --  (p2);
      \draw (p1)  --  (p3);
        \draw (p1)  --  (p4);
     \draw (p1)  --  (p5);
        \draw (p2)  --  (p3);
      \draw (p4)  --  (p5);
        \draw (p3)  --  (p7);
               \draw (p5)  --  (p7);
              \draw (p2)  --  (p6);
               \draw (p4)  --  (p6);   
                 \draw[red](p6)  --  (p7);
           
    \draw[dashed, thin](0,-3.3)--(0,0.5); 
      
 \node [draw=white, fill=white] (b) at (0,-3.5) {(a)};
        \end{tikzpicture}
        \hspace{1cm}
      \begin{tikzpicture}[very thick,scale=0.7]
\tikzstyle{every node}=[circle, draw=black, fill=white, inner sep=0pt, minimum width=4pt];
   
       \path (-0.5,-0.5) node (p1)  {} ;
    \path (0.5,-0.5) node (p2)  {} ;
    \path (-0.5,-1.5) node (p3)  {} ;
     \path (0.5,-1.5) node (p4)  {} ;
     
       \path (-1,0) node (p11)  {} ;
    \path (1,0) node (p22)  {} ;
    \path (-1.5,-2.5) node (p33)  {} ;
     \path (1.5,-2.5) node (p44)  {} ;
     
               \draw[red] (p1)  --  (p2);
      \draw (p1)  --  (p3);
        \draw (p2)  --  (p4);
     \draw[red] (p3)  --  (p4);
      
        \draw[red] (p11)  --  (p22);
      \draw (p11)  --  (p33);
        \draw (p22)  --  (p44);
               \draw[red] (p33)  --  (p44);
            
              \draw (p1)  --  (p11);
               \draw (p2)  --  (p22);   
                 \draw (p3)  --  (p33);
                  \draw (p4)  --  (p44);
    \draw[dashed, thin](0,-3.3)--(0,0.5); 
      
 \node [draw=white, fill=white] (b) at (0,-3.5) {(b)};
        \end{tikzpicture}
                \hspace{1cm}
      \begin{tikzpicture}[very thick,scale=0.7]
\tikzstyle{every node}=[circle, draw=black, fill=white, inner sep=0pt, minimum width=4pt];
   
       \path (-1.5,0) node (p1)  {} ;
    \path (1.5,0) node (p2)  {} ;
    \path (-1.5,-2.5) node (p3)  {} ;
     \path (1.5,-2.5) node (p4)  {} ;
            \path (-0.8,-0.9) node (p5)  {} ;
    \path (0.8,-1.6) node (p6)  {} ;

               \draw (p1)  --  (p5);
      \draw (p5)  --  (p3);
        \draw (p3)  --  (p1);
    
               \draw (p2)  --  (p4);
      \draw (p4)  --  (p6);
        \draw (p6)  --  (p2);
        
              \draw (p1)  --  (p2);
      \draw (p3)  --  (p4);
        \draw[red] (p6)  --  (p5);
          
 \node [draw=white, fill=white] (b) at (0,-3.5) {(c)};
        \end{tikzpicture}
\end{center}
\vspace{-0.6cm} \caption{Symmetric planar frameworks in $\mathbb{R}^2$: (a) and (b) have reflection symmetry and (c) has half-turn symmetry. Bars that are unshifted by the respective reflection or half-turn are shown in red.} \label{fig:symfw}
\end{figure}
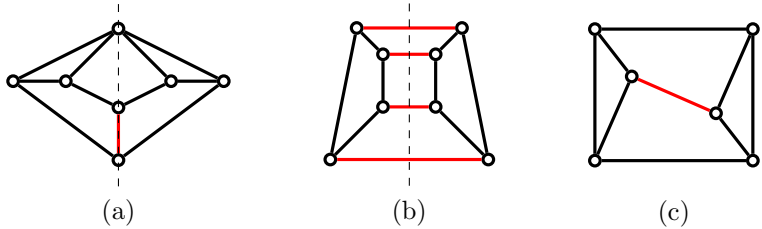

This paper provides methods for designing planar frameworks (or form diagrams) that have additional states of self-stress that cannot be detected with the standard Maxwell count. The nature of these states of self-stress is also investigated with an emphasis on designing fully-symmetric and anti-symmetric states of self-stress.  It is shown that the Fowler-Guest counting rule may be used to increase the number of detected self-stresses and mechanisms of certain symmetry types in either statically determinate or indeterminate frameworks by simply placing a suitable amount of structural members so that they are \emph{unshifted} by the symmetry operations of the framework (see Figure~\ref{fig:symfw}). Due to their simplicity the derived formulas provide a powerful and efficient tool for the design of frameworks with some prespecified structural rigidity properties. As described in Section \ref{sec:formulas}, mirror symmetry plays a larger role than rotational symmetry.  

There are further methods -- beyond symmetry -- that can be used to create additional states of self-stress. These include subdivision methods and tools from projective geometry, such as `pure conditions' \cite{WW,Nixon}, to name but a few. Some of these methods are discussed in Section~\ref{sec:nonsym}  by means of some basic examples. Note that the force-density method \cite{Schek} always produces one state of self-stress but cannot be used to produce any more.

The paper is organised as follows. Section~\ref{sec:prelim} provides a summary of the mathematical background for the symmetry-extended Fowler-Guest counting rule. The formulas for the states of self-stress and mechanisms that can be obtained from this counting rule -- along with a discussion of the key observations arising from these formulas -- are established in Section~\ref{sec:formulas}. These analyses can easily be extended to pinned frameworks as shown in Section~\ref{sec:pinned}. Further methods for creating states of self-stress without using symmetry are discussed in Section~\ref{sec:nonsym}. Finally, we briefly describe some avenues for future research in Section~\ref{sec:future}. 
To demonstrate the applications of our results, examples are given throughout the paper and a detailed discussion of the hypothetical gridshell project shown in Figure~\ref{fig:pic} is given in Section \ref{sec:pinned}.

\section{Preliminaries} \label{sec:prelim}

\subsection{Bar-joint frameworks}

A pin-jointed bar assembly in the plane may be modelled mathematically as a \emph{bar-joint framework}  (or simply {\em framework}) $(G,p)$,
where $G=(V,E)$ is a finite simple graph and  $p:V\to \mathbb{R}^2$ is a map such that $p(i) \neq p(j)$ for all $i,j\in V$. We write each point $p(i)$ as $p_i=(x_i,y_i)$. Each edge of $G$ represents a rigid straight bar and each vertex of $G$ represents a joint or pin that allows rotation in any direction of the plane. We denote $v$ and $e$ to be the number of vertices and edges of $G$, respectively, and throughout this paper we assume that $(G,p)$ is \emph{planar} in the sense that no bars cross each other, and no bar crosses over a joint. Moreover, we  assume that the points of $(G,p)$ affinely span all of the plane.

The \emph{rigidity matrix} $R(G,p)$ of a framework $(G,p)$ is the $e\times 2v$ matrix  
\begin{displaymath} \bordermatrix{& & & & i & & & & j & & & \cr & & & &  & & \vdots & &  & & &
\cr ij & 0 & \ldots &  0 & (p_{i}-p_{j}) & 0 & \ldots & 0 & (p_{j}-p_{i}) &  0 &  \ldots&  0 \cr & & & &  & & \vdots & &  & & &
}
\textrm{,}\end{displaymath}
where, for each edge $ij\in E$ joining the vertices $i$ and $j$, $R(G,p)$ has the row with
$(x_i-x_j)$ and $(y_i-y_j)$ in the two columns associated with $i$,
$(x_j-x_i)$ and $(y_j-y_i)$ in the columns associated with $j$,
and $0$ elsewhere (see, for example, \cite{SW1,Wmatroids}). 

It is well-known that the null-space of $R(G,p)$ is the space of \emph{infinitesimal motions} of $(G,p)$. An infinitesimal motion arising from a rigid body motion in the plane 
 is called a \emph{trivial infinitesimal motion}. The dimension of the space of trivial infinitesimal motions of a framework in the plane is equal to $3$. We will denote the dimension of the space of non-trivial infinitesimal motions, which are often also called \emph{flexes} or \emph{mechanisms}, by $m$. A framework is called \emph{infinitesimally rigid} (or equivalently \emph{statically rigid}) if $m=0$ \cite{Wmatroids}. In structural engineering, an infinitesimally rigid framework is often also called \emph{kinematically determinate} (see \cite{pell} for example).

A \emph{self-stress} of a framework $(G,p)$ is a function  $\omega:E\to \mathbb{R}$ such that for each  vertex $i$ of $G$ the following vector equation holds:
\begin{displaymath}
\sum_{j :ij\in E}\omega(ij)(p_{i}-p_{j})=0 \textrm{.}
\end{displaymath}
In structural engineering, $\omega(ij)(p_{i}-p_{j})$ is called the \emph{axial force} in the bar $ij$, and the stress-coefficient $\omega(ij)$ is called the \emph{force-density} (scalar force divided by the bar length, often written as $T/L$) of the bar $ij$.
 The summation above for vertex $i$ is called the \emph{equilibrium of forces at node $i$}. A self-stress is often also called an \emph{equilibrium stress} as it records  tensions and compressions in the bars balancing at each vertex.

Note that $\omega\in \mathbb{R}^{E}$ is a self-stress if and only if it is a row dependence of $R(G,p)$. 
Equivalently, $\omega\in \mathbb{R}^{E}$ is a self-stress if and only if $R(G,p)^{\top}\omega=0$.
 We will denote the dimension of the space of self-stresses of $(G,p)$ by $s$. A framework with $m=0$ and $s=0$ is called \emph{isostatic}. Isostatic frameworks are minimally infinitesimally rigid and maximally self-stress free.

It follows immediately from the size of the rigidity matrix that a  framework with $e$ edges (or bars) and $v$ vertices (or joints) obeys the Maxwell rule \cite{bibmaxwell} (see also \cite{Cal})
\begin{equation}
    m-s=2v-e-3.
\label{eq:scrule}
\end{equation}

Thus, a necessary condition for a framework to be isostatic is that $e=2v-3$. This condition is not sufficient, however, since a framework may satisfy $e=2v-3$ and $m=s\neq 0$. (See Figure~\ref{fig:mechandstress} for an example.)

\subsection{Block-diagonalisation of the rigidity matrix} \label{sec:symfw}

It was shown in \cite{KG2,KG3} that the rigidity matrix of a  framework $(G,p)$ with point group symmetry $\mathcal{G}$ can be transformed into a block-diagonalised form using methods from group representation theory. In this section we provide the key mathematical background. For the full details, we refer the reader to \cite{owen,Sblock,ST,SW2}.  A \emph{group representation} of $\mathcal{G}$ is a homomorphism from $\mathcal{G}$ to the general linear group of some vector space. The \emph{dimension} of the representation is the dimension of that vector space.

The two key group representations that are needed to obtain the block-decomposition of the rigidity matrix are the `internal' and `external' representation of $(G,p)$ whose corresponding vector spaces are $\mathbb{R}^e$ and $\mathbb{R}^{2v}$ (hence the names `internal' and `external') and which we define below (see also \cite{KG2,KG3,Sblock}). 
Note that each symmetry operation $g\in \mathcal{G}$ of $(G,p)$ induces a permutation of the vertices and bars of $(G,p)$. By a slight abuse of notation, we denote the image of a vertex $i$ or bar $b$ under these permutations by $g(i)$ and $g(b)$, respectively. 

The \emph{internal representation} $P_E:\mathcal{G}\to  GL(\mathbb{R}^e)$ is the permutation representation of the bars of $(G,p)$, that is $P_E(g)= [\delta_{b,g(b')}]_{b,b'}$ for each $g\in \mathcal{G}$, where $\delta$ denotes the Kronecker delta. In other words, the matrix  $P_E(g)$ is the $(0,1)$ matrix which describes how the bars of $(G,p)$ are permuted by $g$.

  Similarly, the \emph{external representation} is defined as $(P_V \otimes T):\mathcal{G}\to GL(\mathbb{R}^{2v})$, where $P_V(g)= [\delta_{i,g(i')}]_{i,i'}$ for each $g\in \mathcal{G}$, $T(g)$ is the matrix in the orthogonal group $O(\mathbb{R}^2)$ representing the isometry $g\in \mathcal{G}$, and $(P_V \otimes T)(g)$ denotes  the Kronecker product of $P_V(g)$ and $T(g)$. In other words, the external representation describes how the vertices are being permuted and how the coordinate system for each vertex is affected by each symmetry operation $g\in\mathcal{G}$.

For a framework  $(G,p)$ with point group symmetry $\mathcal{G}$ we have the following basic intertwining property \cite{Sblock,ST}: $$P_E^{-1}(g) R(G,p) (P_V \otimes T)(g) \quad \textrm{ for all } g\in\mathcal{G}.$$ By Schur's lemma \cite{liebeck,serre}, this implies that the rigidity matrix $R(G,p)$ can be block-decomposed by choosing suitable symmetry-adapted bases for  $\mathbb{R}^{e}$ and $\mathbb{R}^{2v}$. 
More precisely,  if $\rho_1,\ldots, \rho_r$ are the irreducible representations of $\mathcal{G}$, then  the rigidity matrix of $(G,p)$ can be put into the following block form
\begin{equation*}
\label{rigblocks}
A^{\top}R(G,p)B:=\widetilde{R}(G,p)
=\left(\begin{array}{ccc}\widetilde{R}_{1}(G,p) & & \mathbf{0}\\ & \ddots & \\\mathbf{0} &  &
\widetilde{R}_{r}(G,p) \end{array}\right)\textrm{,}
\end{equation*}
where the submatrix block $\widetilde{R}_{i}(G,p)$ corresponds to the irreducible representation $\rho_i$ of $\mathcal{G}$, and $A$ and $B$ are the respective matrices of basis transformation from the standard bases of $\mathbb{R}^e$ and $\mathbb{R}^{2v}$ to the symmetry-adapted bases.

This block-decomposition of the rigidity matrix corresponds to a decomposition  $\mathbb{R}^{e}=X_1 \oplus \cdots \oplus X_r$ of the space $\mathbb{R}^{e}$ into a direct sum of $P_E$-invariant subspaces $X_i$, and a decomposition   $\mathbb{R}^{2v}=Y_1 \oplus \cdots \oplus Y_r$   of the space $\mathbb{R}^{2v}$ into a direct sum of $(P_V \otimes T)$-invariant subspaces $Y_i$, where for a group representation $\Phi:\mathcal{G}\to GL(\mathbb{R}^n)$, a subspace  $U\subseteq \mathbb{R}^n$ is called \emph{$\Phi$-invariant} if $\Phi(g)(U)\subseteq U$ for all $g\in \mathcal{G}$. The spaces $X_i$ and $Y_i$ are associated with $\rho_i$ and the submatrix  $\widetilde{R}_{i}(G,p)$ is of size $\textrm{dim}(X_i)\times \textrm{dim} (Y_i)$. We refer the reader to \cite{Sblock} for the  full mathematical details.

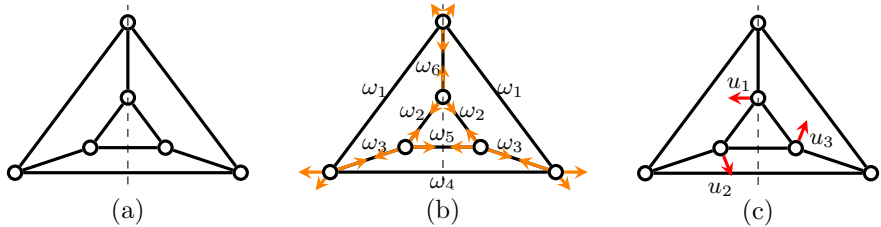
\begin{figure}[htp]
\begin{center}
   \begin{tikzpicture}[very thick,scale=0.5]
\tikzstyle{every node}=[circle, draw=black, fill=white, inner sep=0pt, minimum width=5pt];
        \path (-1,0) node (p1) {} ;
        \path (1,0) node (p2) {} ;
        \path (0,1.33333) node (p3) {} ;
                
         \path (-3,-0.66666) node (p11) {} ;
        \path (3,-0.66666) node (p22) {} ;
        \path (0,3.333333) node (p33) {} ;
           
        \draw (p1)  --  (p2);
         \draw (p1)  --  (p3);
        \draw (p3)  --  (p2);
        \draw (p11)  --  (p22);
         \draw (p11)  --  (p33);
        \draw (p33)  --  (p22);
         \draw (p1)  --  (p11);
         \draw (p2)  --  (p22);
        \draw (p3)  --  (p33);
        
            \draw[dashed, thin](0,-1.4)--(0,3.8);

         \node [draw=white, fill=white] (a) at (0,-1.7) {(a)};
        
                  \end{tikzpicture}
\hspace{0.4cm}
      \begin{tikzpicture}[very thick,scale=0.5]      
\tikzstyle{every node}=[circle, draw=black, fill=white, inner sep=0pt, minimum width=5pt];
        \path (-1,0) node (p1) {} ;
        \path (1,0) node (p2) {} ;
        \path (0,1.33333) node (p3) {} ;
                
         \path (-3,-0.66666) node (p11) {} ;
        \path (3,-0.66666) node (p22) {} ;
        \path (0,3.333333) node (p33) {} ;
               
               \node [draw=white, fill=white] (a) at (-1.8,0.05) {\small$\omega_3$};
          \node [draw=white, fill=white] (a) at (1.8,0.05) {\small$\omega_3$};
            \node [draw=white, fill=white] (a) at (-1.8,1.5) {\small$\omega_1$};
         \node [draw=white, fill=white] (a) at (1.8,1.5) {\small$\omega_1$};
         \node [draw=white, fill=white] (a) at (-0.8,0.9) {\small$\omega_2$};
          \node [draw=white, fill=white] (a) at (0.8,0.9) {\small$\omega_2$};  
               \node [draw=white, fill=white] (a) at (0,-0.98) {\small$\omega_4$};
          \node [draw=white, fill=white] (a) at (-0.37,2.1) {\small$\omega_6$};  
         \node [draw=white, fill=white] (a) at (0,0.3) {\small$\omega_5$};  
               
        \draw (p1)  --  (p2);
         \draw (p1)  --  (p3);
        \draw (p3)  --  (p2);
        \draw (p11)  --  (p22);
         \draw (p11)  --  (p33);
        \draw (p33)  --  (p22);
         \draw (p1)  --  (p11);
         \draw (p2)  --  (p22);
        \draw (p3)  --  (p33);

         \draw[orange,->,>=stealth] (p11)  --  (-3.85,-0.66666);
        \draw[orange,->,>=stealth] (p22)  --  (3.85,-0.6666);
          \draw[orange,->,>=stealth] (p1)  --  (-0.18,0);
        \draw[orange,->,>=stealth] (p2)  --  (0.18,0);
            \draw[orange,->,>=stealth] (p3)  --  (0,2.2);
        \draw[orange,->,>=stealth] (p33)  --  (0,2.46666);
          \draw[orange,->,>=stealth] (p11)  --  (-3.39,-1.15);
        \draw[orange,->,>=stealth] (p33)  --  (0.32,3.82);
         \draw[orange,->,>=stealth] (p22)  --  (3.39,-1.15);
        \draw[orange,->,>=stealth] (p33)  --  (-0.32,3.82);
        
        \draw[orange,->,>=stealth] (p11)  --  (-2,-0.33);
        \draw[orange,->,>=stealth] (p22)  --  (2,-0.33);
         \draw[orange,->,>=stealth] (p1)  --  (-1.85,-0.29);
        \draw[orange,->,>=stealth] (p2)  --  (1.85,-0.29);
        
         \draw[orange,->,>=stealth] (p1)  --  (-0.61,0.5);
        \draw[orange,->,>=stealth] (p2)  --  (0.61,0.5);
        \draw[orange,->,>=stealth] (p3)  --  (-0.35,0.86);
        \draw[orange,->,>=stealth] (p3)  --  (0.35,0.86);
        
            \draw[dashed, thin](0,-1.4)--(0,3.8);

         \node [draw=white, fill=white] (a) at (0,-1.7) {(b)};
        
                  \end{tikzpicture}
\hspace{0.4cm}
  \begin{tikzpicture}[very thick,scale=0.5]
\tikzstyle{every node}=[circle, draw=black, fill=white, inner sep=0pt, minimum width=5pt];
        \path (-1,0) node (p1) {} ;
        \path (1,0) node (p2) {} ;
        \path (0,1.33333) node (p3) {} ;
        
           \node [draw=white, fill=white] (a) at (-0.5,1.7) {\small $u_1$};
          \node [draw=white, fill=white] (a) at (-1,-1.1) {\small$u_2$};  
         \node [draw=white, fill=white] (a) at (1.7,0.2) {\small$u_3$};
               
         \path (-3,-0.66666) node (p11) {} ;
        \path (3,-0.66666) node (p22) {} ;
        \path (0,3.333333) node (p33) {} ;

        \draw (p1)  --  (p2);
         \draw (p1)  --  (p3);
        \draw (p3)  --  (p2);
        \draw (p11)  --  (p22);
         \draw (p11)  --  (p33);
        \draw (p33)  --  (p22);
         \draw (p1)  --  (p11);
         \draw (p2)  --  (p22);
        \draw (p3)  --  (p33);

         \draw[red,->,>=stealth] (p1)  --  (-0.7,-0.75);
        \draw[red,->,>=stealth] (p2)  --  (1.3,0.75);
        \draw[red,->,>=stealth] (p3)  --  (-0.8,1.33333);
        
         \draw[dashed, thin](0,-1.4)--(0,3.8); 
        
            \node [draw=white, fill=white] (a) at (0,-1.7) {(c)};
      \end{tikzpicture}
                \end{center}
\caption{A framework in $\mathbb{R}^2$ with reflection symmetry $\mathcal{C}_s=\{E,\sigma\}$, where $E$ is the identity operation and $\sigma$ is the reflection (a). The framework has a $\rho_1$-symmetric (or `fully-symmetric') self-stress (b) with stress-coefficients $\omega_i$ (the $\omega_i$ are preserved by $\sigma$) and a $\rho_2$-symmetric (or `anti-symmetric') mechanism (c) with velocities $u_i$ (the $u_i$ are reversed by $\sigma$), where $\rho_1(E)=\rho_1(\sigma)=1$ and $\rho_2(E)=1$ and $\rho_2(\sigma)=-1$. 
 In the Mulliken notation \cite{altherz,atkchil} the characters of $\rho_1$ and $\rho_2$ are denoted by $A'=(1,1)$ and $A''=(1,-1)$, respectively.}
\label{fig:mechandstress}
\end{figure}

A vector in $\mathbb{R}^e$ is called \emph{$\rho_i$-symmetric} if it lies in the $P_E$-invariant subspace $X_i$ of $\mathbb{R}^e$. Similarly, a vector in $\mathbb{R}^{2v}$ is called \emph{$\rho_i$-symmetric} if it lies in the $(P_V \otimes T)$-invariant subspace $Y_i$ of $\mathbb{R}^{2v}$. 
See Figure~\ref{fig:mechandstress} for an example.

The space of trivial infinitesimal motions  can be written as the direct sum of the space of translations $\mathcal{T}$ and the space of rotations $\mathcal{R}$, each of which is also a $(P_V \otimes T)$-invariant subspace \cite{Sblock}. Thus, we also have the direct sum decompositions  $\mathcal{T}=T_1 \oplus \cdots \oplus T_r$ and $\mathcal{R}=R_1 \oplus \cdots \oplus R_r$ into $(P_V \otimes T)$-invariant  subspaces $T_i$ and $R_i$, respectively. It follows that for each $i=1,\ldots, r$, we obtain the necessary condition $$\textrm{dim}(X_i)=\textrm{dim}(Y_i)-(\textrm{dim}(T_i)+\textrm{dim}(R_i))$$ for a framework with point group symmetry $\mathcal{G}$ to be isostatic. Using basic results from  character theory, these conditions can be written in a more succinct form as follows \cite{Sblock,owen}:
\begin{equation}
    \Gamma(e)=(\Gamma(v)\times \Gamma_{\mathrm{T}})-(\Gamma_{\mathrm{T}}+\Gamma_{\mathrm{R}}).
\label{eq:bj1a}\end{equation} 

In the terminology of mathematical group theory, each $\Gamma$ in this equation is the character of a group representation of the point group $\mathcal{G}$ of the framework. The \emph{character} of a group representation $\Phi:\mathcal{G}\to \mathbb{R}^n$ associates to each group element of $\mathcal{G}$ the trace of the corresponding matrix (which is independent of the choice of basis for $\mathbb{R}^n$). So for a fixed order of the group elements, the character may be considered as a  $|\mathcal{G}|$-dimensional vector. It is well known that the trace is a class function, so  the entry of the character is the same for each element in the same conjugacy class of the group \cite{liebeck,serre}.
(Note that, confusingly, in applied group theory, a character is usually called a representation, and the trace is called the character, but we will use the mathematical terminology introduced above instead.)
For the point groups in the plane and the characters of their irreducible representations, we will use the standard Schoenflies and Mulliken notations, respectively \cite{altherz,atkchil}.

In equation (\ref{eq:bj1a}), $\Gamma(v)$ and $\Gamma(e)$ are the characters of the permutation representations $P_V$ and $P_E$ of the vertices and edges of $(G,p)$, respectively. That is, the entry of the character $\Gamma(v)$ (or $\Gamma(e)$) corresponding to a group element $g\in \mathcal{G}$ is equal to the number of vertices (edges, respectively) of $(G,p)$ that remain unshifted by the symmetry operation $g$ (since only unshifted structural components contribute a $1$ to the diagonal of the corresponding permutation matrix). See Figure~\ref{fig:symfw} for examples of bars that are unshifted by a reflection or half-turn. In addition, $\Gamma_{\mathrm{T}}$ and $\Gamma_{\mathrm{R}}$ are the  characters of the sub-representation of the external representation  $(P_V \otimes T)$ of $\mathcal{G}$ restricted to the space of translations $\mathcal{T}$ and the space of rotations $\mathcal{R}$, respectively. Note  that $\Gamma(v)\times \Gamma_T= \Gamma(P_V\otimes T)$. 

All the characters in (\ref{eq:bj1}) can be computed by standard manipulations of the character table of the group $\mathcal{G}$ \cite{altherz,atkchil}. See also Table~\ref{tab:2D}.

\subsection{The symmetry-extended Maxwell rule} \label{sec:symfwmaxwell}

The character $\Gamma(\Phi)$ of a group representation $\Phi$ of $\mathcal{G}$ can always be written uniquely as a linear combination of the characters of the irreducible representations $\Gamma(\rho_1),\ldots, \Gamma(\rho_r)$ of $\mathcal{G}$ \cite{liebeck, serre}. It is a standard result in character theory that  the coefficient $\alpha_j$ of each $\Gamma(\rho_j)$ in this linear combination is a non-negative integer and can be found via the following simple formula \cite{liebeck,serre}: 
\begin{equation}\alpha_j=\frac{1}{\|\Gamma(\rho_j)\|^2}\langle \Gamma(\Phi), \Gamma(\rho_j)\rangle,\label{eq:coef}\end{equation}
where $\langle \cdot, \cdot \rangle$ denotes the standard inner product.

Suppose that $(G,p)$ is a framework with point group $\mathcal{G}$ and that $\Gamma(e)=\alpha_1 \Gamma(\rho_1)+\cdots + \alpha_r \Gamma(\rho_r)$ and  $(\Gamma(v)\times \Gamma_T) - (\Gamma_T+\Gamma_R) =\beta_1 \Gamma(\rho_1)+\cdots + \beta_r \Gamma(\rho_r)$, where $\alpha_i,\beta_i\in\mathbb{N}\cup\{0\}$ for all $i=1,\ldots, r$. 
If  $\alpha_i\neq \beta_i$ for some $i$, then it follows from Equation~(\ref{eq:bj1a}) that   $(G,p)$  is not isostatic. 
Moreover, by comparing the coefficients $\alpha_i$ and $\beta_i$ for each $i$, we obtain information about the size of each of the block-matrices $\widetilde{R}_{i}(G,p)$ of the block-decomposed rigidity matrix, which in turn reveals information about the existence of $\rho_i$-symmetric self-stresses or mechanisms. 
So by subtracting $\Gamma(e)$ from $(\Gamma(v)\times \Gamma_{\mathrm{T}})-(\Gamma_{\mathrm{T}}+\Gamma_{\mathrm{R}})$ we obtain the symmetry-extended Maxwell rule, as formulated by Fowler and Guest in \cite{FGsymmax}:
\begin{equation}\label{eq:bj1}\Gamma(m)-\Gamma(s) =(\Gamma(v)\times \Gamma_{\mathrm{T}})-\Gamma(e)-(\Gamma_{\mathrm{T}}+\Gamma_{\mathrm{R}}).\end{equation}
By a slight abuse of terminology, $\Gamma(m)$ and $\Gamma(s)$ are often called the \emph{characters of the mechanisms and states of self-stress} of $(G,p)$, respectively.
If we denote $\gamma_i=\beta_i-\alpha_i$, then $$\Gamma(m)-\Gamma(s)=\sum_{i=1}^{r}\gamma_i\Gamma(\rho_i),$$ where $\gamma_i\in \mathbb{Z}$.  
If $\gamma_i<0$ then we may deduce that $(G,p)$ has a space of $\rho_i$-symmetric self-stresses of dimension at least $-k_i\gamma_i$, where $k_i$ is the dimension of the irreducible representation $\rho_i$ (as defined in the beginning of Section~\ref{sec:symfw}). Similarly, if $\gamma_i>0$ then we may deduce that $(G,p)$ has a space of $\rho_i$-symmetric mechanisms of dimension at least $k_i\gamma_i$.

If there is a mechanism and a self-stress that are both $\rho_i$-symmetric (i.e. they lie in $Y_i$ and $X_i$, respectively), then they cancel in the symmetry-extended count, and can hence not be detected with this count. In particular, we may have $\gamma_i=0$ but $\alpha_i=\beta_i\neq 0$. To find these types of equi-symmetric mechanisms and self-stresses one would have to investigate the null-space and left null-space of the rigidity matrix.

We refer to those  mechanisms and states of self-stress that cannot be detected using the basic Maxwell rule (\ref{eq:scrule}) but are revealed by the symmetry-extended Maxwell rule (\ref{eq:bj1}) as \emph{symmetry-detectable}. Note that for every symmetry-detectable self-stress there exists a symmetry-detectable mechanism and vice versa.


\subsection{Characters for the symmetry-extended Maxwell rule}

The relevant symmetry operations in the plane are: the identity
($E$), rotation by $\phi=2\pi/n$ about a point ($C_n$), and
reflection in a line ($\sigma$). The possible point groups are the infinite set ${\cal C}_n$ and ${\cal C}_{nv}$ for all natural numbers $n$. ${\cal C}_n$ is the cyclic group generated by $C_n$, and ${\cal C}_{nv}$ is the dihedral group generated by a $\{ C_n, \sigma\}$ pair.  The group ${\cal C}_{1v}$ is usually called ${\cal C}_s$.

It was shown in \cite{cfgsw} that the entries of  $\Gamma(m)-\Gamma(s)$ in Equation~(\ref {eq:bj1}) can be computed by keeping track of the  fate of the structural
components of the framework under the various symmetry operations, which in turn depends on
how the joints and bars are placed with respect to the symmetry
elements (i.e., the reflection lines, and the center of rotations, which we may assume to be the origin). The calculations are shown in Table~\ref{tab:2D}, which uses the following notation:
\begin{description}
\item[$v$] is the total number of vertices;
\item[$v_c$] is the number of vertices lying on the centre of
    rotation ($C_{n>2}$ or $C_2$)
    (note that we must have  $v_{c}= 0$ or $1$, since we don't  allow vertices to coincide);
\item[$v_\sigma$] is the number of vertices lying on a given
    mirror line;
\item[$e$] is the total number of edges;
\item[$e_2$] is the number of edges left unshifted by a $C_2$ operation
   (note that if $e_2>1$ then edges cross at the origin, so the framework is non-planar. Note also that $C_n$ with $n>2$     shifts all edges);
\item[$e_\sigma$] is the number of edges unshifted by a
    given reflection (an unshifted edge may lie within, or perpendicular to and centred at the
    mirror line).
\end{description}
Each of the counts above refers to a particular symmetry element, and
any structural component may contribute to one or
more count. For example, a vertex counted in $v_c$ also
contributes to $v_\sigma$ for each mirror line present.
\begin{table}\begin{center}
    \begin{tabular}{l|c c c c}
                        & $\phantom{-}E$       &  $C_{n > 2}$   &
                            $\phantom{-}C_2$   & $\phantom{-}\sigma$  \\ \hline
    $\phantom{=}\Gamma(v)$& $\phantom{-}v$     & $v_c$             &
                            $\phantom{-}v_c$   & $\phantom{-}v_\sigma$\\
    $\phantom{=}\times \Gamma_{T}$&   $\phantom{-}2$     & $2\cos\phi$       &
                            $-2$    & $\phantom{-}0$       \\ \hline
    $=\Gamma(v)\times\Gamma_{T}$
                        & $\phantom{-}2v$      & $2 v_c \cos\phi$  &
                            $-2 v_c$& $\phantom{-}0$       \\
    $\phantom{=}-\Gamma(e)$        & $-e$      & $0$               &
                            $-e_2$  & $-e_\sigma$\\
    $\phantom{=}-(\Gamma_{T}+\Gamma_{R})$
                        & $-3$      & $-2\cos\phi-1$   &
                            $\phantom{-}1$     & $\phantom{-}1$ \\ \hline
    $=\Gamma(m) - \Gamma(s)$
                        & $2v-e-3$  & $2(v_c-1)\cos\phi - 1$ &
                            $-2v_c - e_2 +1$ & $-e_\sigma+1$

    \end{tabular}
    \caption{Calculations of characters for the 2D
        symmetry-extended Maxwell equation (\ref{eq:bj1}). Note that the entries in $\Gamma(m) - \Gamma(s)$ may be non-integers.}
    \label{tab:2D}
    \end{center}
\end{table}

\section{Formulas for creating states of self-stress} \label{sec:formulas}

Throughout this paper it is assumed that $(G,p)$ is a planar framework with point group symmetry $\mathcal{G}$  satisfying $m-s=2v-e-3=k$. The integer $k$ is called the \emph{freedom number} of $(G,p)$. Clearly, if $k<0$ then $(G,p)$ has at least $k$ linearly independent self-stresses, and if $k>0$, then $(G,p)$ has at least $k$ linearly independent mechanisms. For any such frameworks we will now  derive formulas for the number of linearly independent self-stresses (and mechanisms) that can be found with the  symmetry-extended Maxwell rule.

\subsection{Reflection symmetry $\mathcal{C}_s$} \label{sec:mirror}

The reflection group has two irreducible representations, both of which are of dimension $1$. In the Mulliken notation their characters (and the representations themselves) are denoted by $A'$ and $A''$, where $A'=(1,1)$ and $A''=(1,-1)$.

 For a framework with $\mathcal{C}_s$ symmetry satisfying the count $2v-e-3=k$, we obtain from Table~\ref{tab:2D} and Equation (\ref{eq:coef}) that
\begin{equation}\Gamma(m)-\Gamma(s)=(k,-e_\sigma+1)= \frac{k-e_\sigma+1}{2}A' + \frac{k+e_\sigma-1}{2}A''.\label{eq:mirror}\end{equation}
Note that if $k$ is even, then the number of edges, $e$, is odd (for otherwise $2v-e$ is even and hence $k=2v-e-3$ is odd.) Since $e$ is odd, $e_\sigma$ is also odd, because each shifted bar has a mirror copy, so that the number of shifted bars is even. Similarly, if $k$ is odd then $e_\sigma$ is even.

\begin{figure}[htp]
\begin{center}

\begin{tikzpicture}[very thick,scale=0.8]
\tikzstyle{every node}=[circle, draw=black, fill=white, inner sep=0pt, minimum width=4pt];
   
       \path (0,0) node (p1)  {} ;
    \path (0,1) node (p2)  {} ;
    \path (0,-1) node (p3)  {} ;
     \path (-1,0) node (p4)  {} ;
       \path (1,0) node (p5)  {} ;
    \path (-1,0.8) node (p6)  {} ;
    \path (1,0.8) node (p7)  {} ;
     \path (-1.2,-1.1) node (p8)  {} ;
     \path (1.2,-1.1) node (p9)  {} ;
        \path (-0.5,0.5) node (p10)  {} ;
    \path (0.5,0.5) node (p11)  {} ;
     \path (-0.5,-0.5) node (p12)  {} ;
     \path (0.5,-0.5) node (p13)  {} ;

  \draw (p6)  --  (p2);
        \draw (p2)  --  (p7);
      \draw (p6)  --  (p10);
        \draw (p2)  --  (p10);
     \draw (p11)  --  (p2);
        \draw (p11)  --  (p7);
      \draw (p6)  --  (p4);
        \draw (p5)  --  (p7);
     
       \draw (p10)  --  (p1);
        \draw (p1)  --  (p11);
      \draw (p4)  --  (p10);
        \draw (p1)  --  (p10);
     \draw (p11)  --  (p5);
        \draw (p12)  --  (p1);
      \draw (p13)  --  (p1);
        \draw (p4)  --  (p12);
     
         \draw (p13)  --  (p5);
        \draw (p4)  --  (p8);
      \draw (p5)  --  (p9);
        \draw (p8)  --  (p12);
     \draw (p9)  --  (p13);
        \draw (p12)  --  (p3);
      \draw (p13)  --  (p3);
        \draw (p3)  --  (p8);
               \draw (p3)  --  (p9);
    \draw[dashed, thin](0,-1.9)--(0,1.9); 
      \node[rectangle, draw=white](a) at (0,-2.2){(a)};
\end{tikzpicture}
\hspace{0.6cm}
\begin{tikzpicture}[very thick,scale=0.8]
\tikzstyle{every node}=[circle, draw=black, fill=white, inner sep=0pt, minimum width=4pt];
   
       \path (0,0) node (p1)  {} ;
    \path (0,0.8) node (p2)  {} ;
    \path (0,1.2) node (p3)  {} ;
     \path (0,-0.5) node (p4)  {} ;
       \path (0,-1.2) node (p5)  {} ;
    \path (-0.6,0.8) node (p6)  {} ;
    \path (0.6,0.8) node (p7)  {} ;
     \path (-0.9,0) node (p8)  {} ;
     \path (0.9,0) node (p9)  {} ;
        \path (-0.8,-0.5) node (p10)  {} ;
    \path (0.8,-0.5) node (p11)  {} ;
     \path (-1.6,0.2) node (p12)  {} ;
     \path (1.6,0.2) node (p13)  {} ;

  \draw (p1)  --  (p2);
        \draw (p2)  --  (p3);
      \draw (p1)  --  (p4);
        \draw (p4)  --  (p5);
     \draw (p1)  --  (p8);
        \draw (p1)  --  (p9);
      \draw (p6)  --  (p2);
        \draw (p2)  --  (p7);
       \draw (p4)  --  (p10);
        \draw (p4)  --  (p11);

  \draw (p3)  --  (p6);
        \draw (p7)  --  (p3);
      \draw (p6)  --  (p8);
        \draw (p8)  --  (p10);
     \draw (p5)  --  (p10);
        \draw (p5)  --  (p11);
      \draw (p9)  --  (p7);
        \draw (p9)  --  (p11);
       \draw (p12)  --  (p10);
        \draw (p12)  --  (p8);
         \draw (p12)  --  (p6);
        \draw (p13)  --  (p9);
         \draw (p13)  --  (p11);
        \draw (p13)  --  (p7);
     \draw[dashed, thin](0,-1.9)--(0,1.9); 
      \node[rectangle, draw=white](a) at (0,-2.2){(b)};
\end{tikzpicture}
\hspace{0.6cm}
\begin{tikzpicture}[very thick,scale=0.8]
\tikzstyle{every node}=[circle, draw=black, fill=white, inner sep=0pt, minimum width=4pt];
   
       \path (-0.4,0.5) node (p1)  {} ;
    \path (0.4,0.5) node (p2)  {} ;
    \path (-0.7,-0.5) node (p3)  {} ;
     \path (0.7,-0.5) node (p4)  {} ;
       \path (-0.8,1) node (p5)  {} ;
    \path (0.8,1) node (p6)  {} ;
    \path (-0.8,-1) node (p7)  {} ;
     \path (0.8,-1) node (p8)  {} ;
     \path (-1.5,0.8) node (p9)  {} ;
        \path (1.5,0.8) node (p10)  {} ;
    \path (-1.5,-0.8) node (p11)  {} ;
     \path (1.5,-0.8) node (p12)  {} ;

  \draw (p1)  --  (p2);
        \draw (p2)  --  (p4);
      \draw (p4)  --  (p3);
        \draw (p3)  --  (p1);
     \draw (p1)  --  (p2);
        \draw (p1)  --  (p5);
           \draw (p1)  --  (p9);
        \draw (p2)  --  (p6);
        \draw (p2)  --  (p10);
       \draw (p3)  --  (p11);
        \draw (p3)  --  (p7);
       \draw (p4)  --  (p8);
        \draw (p4)  --  (p12);
           \draw (p5)  --  (p6);
        \draw (p6)  --  (p10);
        \draw (p12)  --  (p8);
        \draw (p12)  --  (p10);
        \draw (p8)  --  (p7);
        \draw (p7)  --  (p11);
        \draw (p11)  --  (p9);
        \draw (p9)  --  (p5);
    \draw[dashed, thin](0,-1.9)--(0,1.9); 
      \node[rectangle, draw=white](a) at (0,-2.2){(c)};
\end{tikzpicture}
\end{center}
\vspace{-0.6cm} \caption{Reflection-symmetric frameworks  with an anti-symmetric self-stress (a) and fully-symmetric self-stresses (b), (c). Note that (b) and (c) have four bars that are unshifted by the reflection, whereas (a) has none. See Example~\ref{ex:cs} for a more detailed discussion.} \label{fig:anti-sym}
\end{figure}
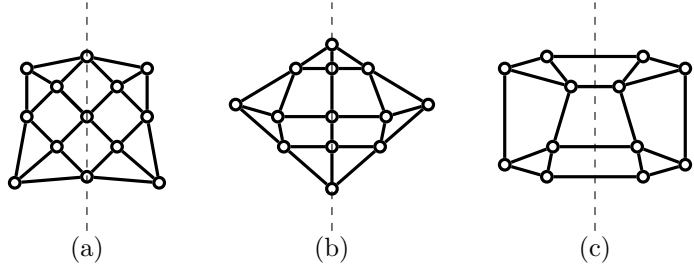

Some observations arising from Equation~(\ref{eq:mirror}) are:
\begin{enumerate}
\item Suppose $k\leq 0$. Then the standard Maxwell rule (\ref{eq:scrule}) tells us that the framework has at least $-k$ linearly independent self-stresses. Note that the coefficients of $A'$ and $A''$ in Equation~(\ref{eq:mirror}) are integers and add up to $k$, and the coefficient of $A'$ is non-positive for any value of $e_\sigma$, since $e_\sigma\geq 0$. If $e_\sigma\leq -k+1$, then the coefficient of $A''$ is also non-positive, and hence we still detect only $-k$ independent self-stresses. However, we may deduce from Equation~(\ref{eq:mirror}) that in this case we have $\frac{-k+e_\sigma-1}{2}$ independent self-stresses that are $A'$-symmetric, and $\frac{-k-e_\sigma+1}{2}$ independent self-stresses that are $A''$-symmetric. (See Figure~\ref{fig:anti-sym}(a) for an example.) 

By definition of the internal representation, the $A'$-symmetric self-stresses are `fully-symmetric' in the sense that mirror images of bars have the same stress-coefficients (recall Figure~\ref{fig:mechandstress}). The $A''$-symmetric self-stresses are `anti-symmetric' in the sense that if an edge has stress-coefficient $\omega$, then its symmetric copy under the reflection has stress-coefficient $-\omega$. (So in particular, the stress-coefficient of any edge that is unshifted by the mirror is zero.)

\item Suppose again that $k\leq 0$. The larger we make $e_\sigma$ while keeping $k$ fixed, the more anti-symmetric self-stresses are switched to fully-symmetric self-stresses. When $e_\sigma=-k+1$ then all $-k$ detected self-stresses are fully-symmetric. If we  increase $e_\sigma$ further so that $e_\sigma\geq -k+3$, then the coefficient  of $A''$ becomes positive and hence we obtain symmetry-detectable $A''$-symmetric mechanisms and, simultaneously, symmetry-detectable $A'$-symmetric self-stresses. So in this case we detect $\frac{-k+e_\sigma-1}{2}>-k$ self-stresses. (See Figure~\ref{fig:anti-sym}(b) for an example.) The more bars are positioned so that they are unshifted by the mirror, while keeping $k$ fixed, the more symmetry-detectable fully-symmetric self-stresses are obtained.

\item In the special case of $k=0$ there are no symmetry-detectable self-stresses or mechanisms if $e_\sigma=1$. In fact, in this case the framework is isostatic for any `generic' positions of the vertices, as shown in \cite{BS4}.  If $e_\sigma\geq 3$, then we obtain $\frac{e_\sigma-1}{2}$ symmetry-detectable fully-symmetric self-stresses.  

\item Suppose $k>0$. Then the coefficient of $A''$ is always non-negative. If $e_\sigma\leq k+1$ then we only find the $k$ mechanisms that were already predicted by the standard Maxwell rule (\ref{eq:scrule}). However, we obtain some valuable information about their symmetry properties.  If $e_\sigma\geq k+3$, then we obtain symmetry-detectable self-stresses, all of which are fully-symmetric. (See Figure~\ref{fig:anti-sym}(c) for an example.) 
\end{enumerate}

In summary, we increase the number of fully-symmetric self-stresses for a fixed $k$ by increasing $e_\sigma$. We increase the number of anti-symmetric self-stresses by decreasing $e_\sigma$.  Note, however, that $e_{\sigma}$ can never be negative.
 
\begin{example} \label{ex:cs}
Figure~\ref{fig:anti-sym} shows three examples of  frameworks with $\mathcal{C}_s$ symmetry. 
The framework in (a)  has $e=2v-2=24$, so $k=-1$, and $e_\sigma=0$. Thus, by Equation~(\ref{eq:mirror}), we have $\Gamma(m)-\Gamma(s)= -A''$. So we only find the self-stress which is guaranteed to exist by the $k=-1$ count, but we see that it is anti-symmetric. 

The framework in (b) has the same underlying graph as the one in (a) and has $k=-1$ and  $e_\sigma=4$. Thus, by Equation~(\ref{eq:mirror}), we have $\Gamma(m)-\Gamma(s)= -2A'+A''$. Since each negative coefficient indicates self-stresses and each positive coefficient indicates mechanisms,  we deduce that the framework has two independent fully-symmetric self-stresses -- one of which is symmetry-detectable -- and one symmetry-detectable anti-symmetric mechanism.

Finally, the framework in (c)  has $e=2v-4=20$, so $k=1$, and $e_\sigma=4$. Thus, by Equation~(\ref{eq:mirror}), we have $\Gamma(m)-\Gamma(s)= -A'+2A''$. It follows that the framework has a symmetry-detectable fully-symmetric self-stress and two anti-symmetric mechanisms, one of which is also symmetry-detectable.
\end{example}

\begin{remark}
As shown in Section~\ref{sec:symfwmaxwell}, the character counts describe the dimensions of the block matrices in the block-decomposed rigidity matrix $\widetilde{R}(G,p)$. There are some standard methods and algorithms for finding the symmetry-adapted bases that give this block-demposition of $\widetilde{R}(G,p)$ (see, for example, \cite{faessler,weeny}). From the specific entries of the block matrices $\widetilde{R}_i(G,p)$, we may then compute their kernels and co-kernels and hence obtain the complete information about the mechanisms and self-stresses of $(G,p)$ and their symmetry types. Recent work has also established `orbit matrices' that are equivalent to the block-matrices and whose entries can be written down directly from the coordinates of the points \cite{SWorbit,ST}. This reduces the computational effort in analysing these matrices. However,  analyses of the kernels or co-kernels of the block-matrices  often do not help the designer in obtaining realisations of graphs with additional states of self-stress.  
\end{remark}

\begin{figure}[htp]
\begin{center}
\begin{tikzpicture}[very thick,scale=1]
\tikzstyle{every node}=[circle, draw=black, fill=white, inner sep=0pt, minimum width=4pt];
   
       \path (-0.7,0) node (p1)  {} ;
    \path (0.7,0) node (p2)  {} ;
           \draw (p1)  --  (p2);
             \draw [->,>=stealth,red] (p1) --(-1.2,0.5);   
                \draw [->,>=stealth,red] (p2) --(1.2,0.5);    
    \draw[dashed, thin](0,-1.3)--(0,1.3); 
      
 \node [draw=white, fill=white] (b) at (0,-1.5) {(a)};
        \end{tikzpicture}
         \hspace{0.7cm}
      \begin{tikzpicture}[very thick,scale=1]
\tikzstyle{every node}=[circle, draw=black, fill=white, inner sep=0pt, minimum width=4pt];
   
       \path (0,-0.5) node (p1)  {} ;
    \path (0,0.5) node (p2)  {} ;
           \draw (p1)  --  (p2);
     \draw [->,>=stealth,red] (p2) --(0,1);   
                \draw [->,>=stealth,red] (p1) --(0,-1);   
    \draw[dashed, thin](0,-1.3)--(0,1.3); 
      
 \node [draw=white, fill=white] (b) at (0,-1.5) {(b)};
        \end{tikzpicture}
                \hspace{0.7cm}
      \begin{tikzpicture}[very thick,scale=1]
\tikzstyle{every node}=[circle, draw=black, fill=white, inner sep=0pt, minimum width=4pt];
   
       \path (0,-0.5) node (p1)  {} ;
    \path (0,0.5) node (p2)  {} ;
           \draw (p1)  --  (p2);
     \draw [->,>=stealth,red] (p2) --(0,1);   
                \draw [->,>=stealth,red] (p1) --(0,0);   
    \draw[dashed, thin](0,-1.3)--(0,1.3); 
      
 \node [draw=white, fill=white] (b) at (0,-1.5) {(c)};
        \end{tikzpicture}
        \hspace{0.7cm}
        \begin{tikzpicture}[very thick,scale=1]
\tikzstyle{every node}=[circle, draw=black, fill=white, inner sep=0pt, minimum width=4pt];
   
       \path (-0.7,0) node (p1)  {} ;
    \path (0.7,0) node (p2)  {} ;
           \draw (p1)  --  (p2);
             \draw [->,>=stealth,red] (p1) --(-1.2,0.5);   
                \draw [->,>=stealth,red] (p2) --(0.2,-0.5);    
    \draw[dashed, thin](0,-1.3)--(0,1.3); 
      
 \node [draw=white, fill=white] (b) at (0,-1.5) {(d)};
        \end{tikzpicture}
               \hspace{0.7cm}
      \begin{tikzpicture}[very thick,scale=1]
\tikzstyle{every node}=[circle, draw=black, fill=white, inner sep=0pt, minimum width=4pt];
   
       \path (0,-0.5) node (p1)  {} ;
    \path (0,0.5) node (p2)  {} ;
           \draw (p1)  --  (p2);
     \draw [->,>=stealth,red] (p1) --(-0.8,-0.5);   
                \draw [->,>=stealth,red] (p2) --(0.8,0.5);   
    \draw[dashed, thin](0,-1.3)--(0,1.3); 
      
 \node [draw=white, fill=white] (b) at (0,-1.5) {(e)};
        \end{tikzpicture}
\end{center}
\vspace{-0.6cm} \caption{Velocity vectors at the vertices of a bar that is unshifted by a mirror. The velocities in (a), (b) and (c) are fully-symmetric. The ones in (a) and (b) do not form an infinitesimal motion, since their orthogonal projections onto the bar create a non-zero strain on the bar. The velocities in (d) and (e) are anti-symmetric. Note that \emph{any} anti-symmetric velocity assignment will yield an infinitesimal motion of an unshifted bar, and hence such a bar does not impose any constraint when restricting to anti-symmetric velocity assignments.} \label{fig:loops}
\end{figure}
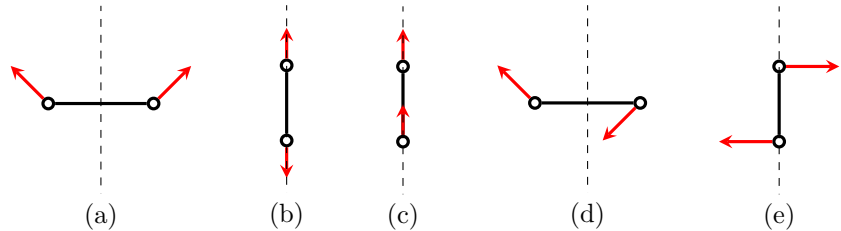

\begin{remark}
 The above observations on numbers of self-stresses and their symmetry types are a consequence of the fact that a bar that is unshifted by a reflection does not constitute any constraint when we restrict to anti-symmetric assignments of velocity vectors. See Figure~\ref{fig:loops} for an illustration. So an unshifted bar always contributes a row to the fully-symmetric block matrix of $\widetilde{R}(G,p)$, and not to the anti-symmetric one. 
 
 Thus, if we start with a fixed freedom number $k$ and increase $e_\sigma$ then we may create additional row dependencies in the fully-symmetric block matrix and remove row dependencies in the anti-symmetric block matrix, but not vice versa. In other words, by increasing $e_\sigma$, we can only switch anti-symmetric self-stresses to fully-symmetric ones. 
\end{remark}

\subsection{Half-turn symmetry $\mathcal{C}_2$}

The half-turn rotational group has two irreducible representations, which are the same as for the reflection group. These representations and their characters are denoted by $A=(1,1)$ and $B=(1,-1)$.

For a framework with $\mathcal{C}_2$ symmetry satisfying the count $2v-e-3=k$, we obtain from Table~\ref{tab:2D} and Equation~(\ref{eq:coef}) that
\begin{equation}\Gamma(m)-\Gamma(s)=(k,-2v_c-e_2+1)= \frac{k-2v_c-e_2+1}{2}A + \frac{k+2v_c+e_2-1}{2}B.\label{eq:halfturn}\end{equation}
Note that if $k$ is even, then $e_2$ is odd. Similarly, if $k$ is odd then $e_2$ is even.

It follows from the definition of a framework that $v_c$ (i.e, the number of vertices of $(G,p)$ positioned at the origin) equals $0$ or $1$. Moreover, by our planarity assumption, if $v_c=0$ then $e_2=0$ or $1$, and if $v_c=1$ then $e_2=0$. Thus, Equation~(\ref{eq:halfturn}) simplifies to
\begin{equation}\label{eq:halfturn1} \Gamma(m)-\Gamma(s)=\begin{cases}
\frac{k+1}{2}A + \frac{k-1}{2}B, &  \textrm{ if } v_c=e_2=0\\
\frac{k}{2}A + \frac{k}{2}B, &  \textrm{ if } v_c=0,\, e_2=1\\
\frac{k-1}{2}A + \frac{k+1}{2}B, &  \textrm{ if } v_c=1\\
\end{cases}\end{equation}
Some observations arising from Equation~(\ref{eq:halfturn1}) are:
\begin{enumerate}
\item If $v_c=0$, then we obtain the same count for $\Gamma(m)-\Gamma(s)$ as we did for reflection symmetry, with $e_\sigma$ being replaced by $e_2$. (This is not surprising, given the transfer results for infinitesimal rigidity between $\mathcal{C}_s$ and $\mathcal{C}_2$ established in \cite{cnsw}.) So the same observations we made for the reflection symmetry also apply to the half-turn symmetry in the case when $v_c=0$. However, note that $e_2$ cannot be larger than $1$ by our planarity assumption, whereas $e_\sigma$ did not have this restriction. So unlike in the reflection symmetry case, we cannot keep increasing the number of  fully-symmetric self-stresses by increasing $e_2$.
\item If $v_c=1$, then $e_2=0$ and hence $k$ is odd. In this case, there are no symmetry-detectable mechanisms or self-stresses, since the coefficients of $A$ and $B$ add up to $k$ and are either both non-positive or both non-negative.  If $k=-1$ then we obtain one fully-symmetric self-stress, but no anti-symmetric self-stress. By taking larger negative $k$ we increase both the number of fully-symmetric and anti-symmetric self-stresses. For any positive $k$ we only find mechanisms.
\end{enumerate}

\subsection{Rotational symmetry $\mathcal{C}_n$, $n\geq 3$}

The group $\mathcal{C}_n$ has $n$ irreducible $1$-dimensional representations whose characters are denoted by $A_t$ for $t=0,\ldots, n-1$. The $j$-th entry of the character $A_t$ is given by $(A_t)_j=\epsilon^{tj}$, where $\epsilon$ denotes the complex root of unity $e^{\frac{2\pi i}{n}}$.

Suppose we are given a framework $(G,p)$ with freedom number $k$ and $\mathcal{C}_n$ symmetry, where $n\geq 3$. Then $e_n=0$ and, by definition of a framework, $v_c$ equals $0$ or $1$. Note that if $n$ is even, the group $\mathcal{C}_n$ contains the half-turn $C_n^{n/2}$. However, if $e_2>0$, then the $\mathcal{C}_n$ symmetry implies that $e_2>1$, contradicting the planarity of $(G,p)$. Thus, $e_n=e_2=0$. 
From Table~\ref{tab:2D} we obtain the following.

 For $v_c=0$ we  have:
\begin{equation*}\label{eq:cn2}   \Gamma(m)-\Gamma(s)=\Big(k,-2\cos \frac{2\pi}{n}-1, \ldots,  -2\cos \pi-1,\ldots, -2\cos \frac{(n-1)2\pi}{n}-1\Big).
\end{equation*}
Note here that the entry $-2\cos \pi-1=1$ appears if and only if $n$ is even.

For $v_c=1$  we  have:
\begin{equation*}\label{eq:cn3}   \Gamma(m)-\Gamma(s)=(k,-1,-1,\ldots, -1).
\end{equation*}

In the case when $v_c=0$, we may write $k$ as $k+3-2\cos0-1$. 
 Similarly, in the case when $v_c=1$ we may write $k$ as $k+1-1$. From Equation~(\ref{eq:coef}) and the standard fact that for $t\neq 0$ we have
$\sum_{j=0}^{n-1} \epsilon^{tj}=0$,
 we then obtain the following expressions for $\Gamma(m)-\Gamma(s)$.

\begin{itemize}
\item For $v_c=0$  we  obtain:
\begin{equation*}\Gamma(m)-\Gamma(s)=\Big(\frac{k+3}{n}-1\Big) A_0 + \sum_{t=1}^{n-1} \Big(\frac{k+3-2\sum_{j=0}^{n-1} \epsilon^{tj}\cos \big(\frac{j2\pi}{n}\big)}{n}\Big) A_t
\end{equation*}
which simplifies (by Proposition~\ref{prop:triv} in the Appendix) to 
\begin{equation}\label{eq:cn11}\Gamma(m)-\Gamma(s)=\Big(\frac{k+3}{n}-1\Big) A_0 + \Big(\frac{k+3}{n}-1\Big)A_1+ \sum_{t=2}^{n-2} \frac{k+3}{n} A_t +  \Big(\frac{k+3}{n}-1\Big)A_{n-1}
\end{equation}

\item For $v_c=1$  we  obtain:
\begin{equation}\label{eq:cn33} \Gamma(m)-\Gamma(s)= \Big(\frac{k+1}{n}-1\Big) A_0 + \sum_{t=1}^{n-1} \Big(\frac{k+1}{n}\Big) A_t
\end{equation}
\end{itemize}

Some observations arising from Equations~(\ref{eq:cn11}) and (\ref{eq:cn33}) are: 
\begin{enumerate}
\item If $v_c=0$, then $n$ must divide $k+3$. 
By Equation~(\ref{eq:cn11}), the symmetry-extended counting rule does not detect any self-stresses or mechanisms in addition to the ones that are detected by the standard Maxwell rule. To see this, note that the sum of the coefficients of the $A_t$ equals $k$ and the coefficients are either all non-positive or all non-negative.  Equation~(\ref{eq:cn11}) shows that in the presence of symmetry,  the self-stresses distribute across the $P_E$-invariant subspaces $X_t$ corresponding to $A_t$ as follows. Let $\ell\geq 0$ and  $k=-\ell n-3$. Then we detect $(\ell+1)$   self-stresses of symmetry $A_0$, $A_1$ and $A_{n-1}$, and $\ell$ $A_t$-symmetric self-stresses  for each $t\neq 0,1,n-1$.

\item If $v_c=1$, then $n$ must divide $k+1$. Again, there are no symmetry-detectable self-stresses or mechanisms. Equation~(\ref{eq:cn33}) shows that in the presence of symmetry,  the self-stresses distribute across the $P_E$-invariant subspaces $X_t$ as follows. Let $\ell\geq 0$ and  $k=-\ell n-1$. Then we detect $(\ell+1)$  fully-symmetric self-stresses, and $\ell$ $A_t$-symmetric self-stresses  for each $t\neq 0$.
\end{enumerate} 

In summary, it turns out that for a framework with rotational symmetry $\mathcal{C}_n$, $n\geq 3$, there are no symmetry-detectable self-stresses or mechanisms. Any self-stresses are distributed equally across the different symmetry types $A_t$, except for an extra self-stress of symmetry $A_0$, $A_1$ and $A_{n-1}$ in the case when $v_c=0$, and an extra self-stress of symmetry $A_0$ in the case when $v_c=1$.

\begin{remark} It was shown in \cite[Lemma 6.7]{ST} that 
 the block-matrices of the block-decomposed rigidity matrix $\widetilde{R}(G,p)$ corresponding to $A_1$ and $A_{n-1}$ have a kernel of dimension at least 1, since we may choose a basis for the space of infinitesimal translations that consists of an $A_1$-symmetric and an $A_{n-1}$-symmetric translation. The trivial infinitesimal rotation is $A_0$-symmetric. (See also \cite{IS,BS3,ST} for combinatorial characterisations of infinitesimally rigid frameworks with $\mathcal{C}_n$ symmetry in the case when $v_c=0$ and $n$ is odd.)

So we may interpret Equation~(\ref{eq:cn11}) as follows: if $v_c=0$ and $e_n=e_2=0$, then each block matrix of $\widetilde{R}(G,p)$ has the same size (or, in other words, each edge orbit under the $\mathcal{C}_n$ symmetry contributes one edge to each of the $n$ blocks, and each vertex orbit contributes 2 columns -- or one vertex --  to each of the $n$ blocks), and the extra self-stresses for the blocks corresponding to $A_0, A_1$ and $A_{n-1}$ appear due to the symmetric decomposition of the trivial motion space.

 In the case when $v_c=1$, we have a special vertex orbit of size 1 (the vertex at the center of rotation), which adds one column to each of the blocks corresponding to $A_1$ and $A_{n-1}$ so that we only obtain an extra self-stress for the block corresponding to $A_0$.
\end{remark}

\subsection{Dihedral symmetry $\mathcal{C}_{2v}$} \label{sec:c2v}

Recall that the group $\mathcal{C}_{2v}$ consists of the identity $E$, two reflections $\sigma_h$ and $\sigma_v$ in perpendicular mirror lines, and the half-turn $C_2$. This point group symmetry appears frequently in engineering designs. The characters of the four irreducible representations of $\mathcal{C}_{2v}$  are shown in Table~\ref{tab:c2v}. 
 \begin{table}[htp]
\begin{center}\begin{tabular}{c||c|c|c|c}
$\mathcal{C}_{2v}$   &   $E$  &  $C_{2}$   & $\sigma_{h}$   &    $\sigma_{v}$\\\hline\hline
$A_{1}$  &    1  &  1 &  1  &  1\\\hline
$A_{2}$  &    1  &  1 &  -1  &  -1\\\hline
$B_{1}$  &    1  &  -1 &  1  &  -1\\\hline
$B_{2}$  &    1  &  -1 &  -1  &  1\\
\end{tabular}\caption{The irreducible characters of $\mathcal{C}_{2v}$.} \label{tab:c2v}\end{center}
\end{table}

For a framework with $\mathcal{C}_{2v}$ symmetry satisfying the count $2v-e-3=k$, we obtain from Table~\ref{tab:2D}  that
$$\Gamma(m)-\Gamma(s)= (k,-2v_c-e_2+1,-e_{\sigma_h}+1,-e_{\sigma_v}+1).$$
Thus, by Equation (\ref{eq:coef}) we obtain the following expressions for $\Gamma(m)-\Gamma(s)$.

\begin{itemize}
\item For $v_c=0$ and $e_2=0$ we  obtain:
{
\medmuskip=0mu
\thinmuskip=0mu
\thickmuskip=0mu
\nulldelimiterspace=0pt
\scriptspace=0pt
\begin{equation}\label{eq:c2vcase1} \hspace{-0.7cm}\resizebox{0.9\hsize}{!}{$\Gamma(m)-\Gamma(s)=\frac{k-e_{\sigma_h}-e_{\sigma_v}+3}{4}A_1 + \frac{k+e_{\sigma_h}+e_{\sigma_v}-1}{4}A_2+ \frac{k-e_{\sigma_h}+e_{\sigma_v}-1}{4}B_1 +\frac{k+e_{\sigma_h}-e_{\sigma_v}-1}{4}B_2  $ }
\end{equation}
}
\item For $v_c=0$ and $e_2=1$ we obtain:{
\medmuskip=0mu
\thinmuskip=0mu
\thickmuskip=0mu
\nulldelimiterspace=0pt
\scriptspace=0pt
\begin{equation}\label{eq:c2vcase2}\hspace{-0.7cm}\resizebox{0.9\hsize}{!}{$\Gamma(m)-\Gamma(s)=\frac{k-e_{\sigma_h}-e_{\sigma_v}+2}{4}A_1 + \frac{k+e_{\sigma_h}+e_{\sigma_v}-2}{4}A_2+ \frac{k-e_{\sigma_h}+e_{\sigma_v}}{4}B_1 +\frac{k+e_{\sigma_h}-e_{\sigma_v}}{4}B_2$}   
\end{equation}
}
\item For $v_c=1$ and $e_2=0$ we obtain:
{
\medmuskip=0mu
\thinmuskip=0mu
\thickmuskip=0mu
\nulldelimiterspace=0pt
\scriptspace=0pt
\begin{equation}\label{eq:c2vcase3} \hspace{-0.7cm}
\resizebox{0.9\hsize}{!}{$\Gamma(m)-\Gamma(s)=\frac{k-e_{\sigma_h}-e_{\sigma_v}+1}{4}A_1 + \frac{k+e_{\sigma_h}+e_{\sigma_v}-3}{4}A_2+ \frac{k-e_{\sigma_h}+e_{\sigma_v}+1}{4}B_1 +\frac{k+e_{\sigma_h}-e_{\sigma_v}+1}{4}B_2$}   \end{equation}
}
\end{itemize}

Note that if $k$ is even, then $e$ is odd and so $e_2=1$ and both $e_{\sigma_h}$ and $e_{\sigma_v}$ are odd. Similarly, if $k$ is odd, then $e$ is even and so $e_2=0$ and both $e_{\sigma_h}$ and $e_{\sigma_v}$ are even. For $\mathcal{C}_{2v}$, the notation  $\sigma_h$ and $\sigma_v$ is used for reflections in a horizontal and vertical mirror line, respectively.

In the following we will assume that $e_{\sigma_h}\geq e_{\sigma_v}$.  Some observations arising from Equations~(\ref{eq:c2vcase1})--(\ref{eq:c2vcase3}) are: 
\begin{enumerate}
\item Suppose $k\leq 0$ and $k$ is even. Then we need to consider Equation~(\ref{eq:c2vcase2}). (The analysis for the case when $k\leq 0$ is odd is  analogous, but we need to consider Equation~(\ref{eq:c2vcase1}) or (\ref{eq:c2vcase3}) depending on whether $v_c=0$ or 1.) We denote the coefficients of $A_i$ by $\alpha_i$, and  the coefficients of $B_i$  by $\beta_i$ for $i=1,2$. Note that $\alpha_1+\alpha_2+\beta_1+\beta_2=k$, and that $\alpha_1\leq 0$ and $\beta_1\leq 0$ for any values of $e_{\sigma_h}$ and $e_{\sigma_v}$. 

If $\alpha_2\leq 0$ and $\beta_2\leq 0$, then there are no symmetry-detectable self-stresses: we only find the $-k$ self-stresses predicted by the standard Maxwell rule. Since $\alpha_2\geq \beta_2$, this happens when $\alpha_2\leq 0$, or $e_{\sigma_h}+e_{\sigma_v}-2\leq -k$. However, in this case we still obtain valuable information about the symmetry types of these self-stresses.

Suppose $\alpha_2> 0$, or equivalently, $e_{\sigma_h}+e_{\sigma_v}-2> -k$. We have $\beta_2\geq 0$ if and only if $e_{\sigma_h}-e_{\sigma_v}\geq -k$.  In this case, we detect $-\alpha_1-\beta_1=\frac{-k+e_{\sigma_h}-1}{2}$ self-stresses (and $\alpha_2+\beta_2$ mechanisms). Since $\beta_2\geq 0$ also implies that $e_{\sigma_h}\geq -k+1$, an analysis of the framework using only the reflection symmetry $\mathcal{C}_s$ with mirror $e_{\sigma_h}$ detects the same number of self-stresses as the $\mathcal{C}_{2v}$ analysis (recall Section~\ref{sec:mirror}). However, an analysis with the larger $\mathcal{C}_{2v}$ group again provides added information regarding the symmetry types of the self-stresses. (See Figure~\ref{fig:2cvstresses}(a) for an example.)

Suppose $\alpha_2> 0$ and $\beta_2< 0$, i.e., $e_{\sigma_h}-e_{\sigma_v}< -k< e_{\sigma_h}+e_{\sigma_v}-2$. Then we detect $-\alpha_1-\beta_1-\beta_2=\frac{-3k+e_{\sigma_h}+e_{\sigma_v}-2}{4}>-k$ self-stresses (and $\alpha_2$ mechanisms).  
In this case the $\mathcal{C}_{2v}$ analysis detects more self-stresses than a $\mathcal{C}_s$ analysis, since $\beta_2<0$ implies that $-\alpha_1-\beta_1-\beta_2>\frac{-k+e_{\sigma_h}-1}{2}$.  Note that there will be at least one fully-symmetric self-stress, as well as at least one $B_1$-symmetric and at least one $B_2$-symmetric self-stress in this case, since $\beta_2<0$ implies $\beta_1<0$ and $\alpha_1<0$. So in particular, for each mirror there will be at least one anti-symmetric self-stress. (See Figure~\ref{fig:2cvstresses}(b) for an example.)

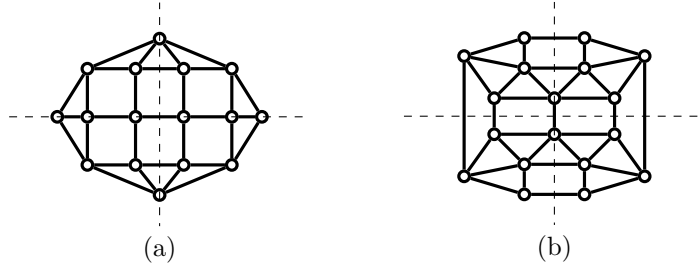
\begin{figure}[htp]
\begin{center}
\begin{tikzpicture}[very thick,scale=0.8]
\tikzstyle{every node}=[circle, draw=black, fill=white, inner sep=0pt, minimum width=4pt];
   
       \path (-0.4,0) node (p1)  {} ;
    \path (0.4,0) node (p1r)  {} ;
            \path (-0.4,0.8) node (p2)  {} ;
    \path (0.4,0.8) node (p2r)  {} ;
       \path (-0.4,-0.8) node (p3)  {} ;
    \path (0.4,-0.8) node (p3r)  {} ;
       \path (0,1.3) node (p4)  {} ;
    \path (0,-1.3) node (p5)  {} ;
       \path (-1.2,0.8) node (p6)  {} ;
    \path (1.2,0.8) node (p6r)  {} ;
       \path (-1.2,0) node (p7)  {} ;
    \path (1.2,0) node (p7r)  {} ;
       \path (-1.2,-0.8) node (p8)  {} ;
    \path (1.2,-0.8) node (p8r)  {} ;
       \path (-1.7,0) node (p9)  {} ;
    \path (1.7,0) node (p9r)  {} ;
     
  \draw (p1)  --  (p2);
        \draw (p1)  --  (p3);
      \draw (p1)  --  (p7);
        \draw (p4)  --  (p2);
     \draw (p6)  --  (p2);
        \draw (p3)  --  (p5);
           \draw (p3)  --  (p8);
      \draw (p5)  --  (p8);
           \draw (p8)  --  (p7);
           \draw (p7)  --  (p6);
           \draw (p4)  --  (p6);
           \draw (p7)  --  (p9);
           \draw (p6)  --  (p9);
                   \draw (p8)  --  (p9);
           
        \draw (p1r)  --  (p2r);
        \draw (p1r)  --  (p3r);
      \draw (p1r)  --  (p7r);
        \draw (p4)  --  (p2r);
     \draw (p6r)  --  (p2r);
        \draw (p3r)  --  (p5);
           \draw (p3r)  --  (p8r);
      \draw (p5)  --  (p8r);
           \draw (p8r)  --  (p7r);
           \draw (p7r)  --  (p6r);
           \draw (p4)  --  (p6r);
           \draw (p7r)  --  (p9r);
           \draw (p6r)  --  (p9r);
                   \draw (p8r)  --  (p9r);
                   
       \draw (p2r)  --  (p2);
           \draw (p1r)  --  (p1);
                   \draw (p3r)  --  (p3);

    \draw[dashed, thin](0,-1.9)--(0,1.9); 
 \draw[dashed, thin](-2.5,0)--(2.5,0); 
 \node [draw=white, fill=white] (b) at (0,-2.2) {(a)};
\end{tikzpicture}
\hspace{1cm}
\begin{tikzpicture}[very thick,scale=0.8]
\tikzstyle{every node}=[circle, draw=black, fill=white, inner sep=0pt, minimum width=4pt];
   
       \path (0,0.3) node (p1)  {} ;
    \path (0,-0.3) node (p2)  {} ;
    \path (-0.5,0.8) node (p3)  {} ;
     \path (0.5,0.8) node (p4)  {} ;
        \path (-0.5,1.3) node (p5)  {} ;
     \path (0.5,1.3) node (p6)  {} ;
             \path (-0.5,-0.8) node (p7)  {} ;
     \path (0.5,-0.8) node (p8)  {} ;
        \path (-0.5,-1.3) node (p9)  {} ;
     \path (0.5,-1.3) node (p10)  {} ;
      
       \path (-1,0.3) node (p11)  {} ;
    \path (-1,-0.3) node (p12)  {} ;
        \path (-1.5,1) node (p15)  {} ;
     \path (-1.5,-1) node (p16)  {} ;
           
       \path (1,0.3) node (p11r)  {} ;
    \path (1,-0.3) node (p12r)  {} ;
        \path (1.5,1) node (p15r)  {} ;
     \path (1.5,-1) node (p16r)  {} ;
     
  \draw (p1)  --  (p2);
        \draw (p1)  --  (p4);
      \draw (p1)  --  (p3);
        \draw (p7)  --  (p2);
     \draw (p8)  --  (p2);
        \draw (p3)  --  (p5);
           \draw (p4)  --  (p6);
      \draw (p7)  --  (p9);
           \draw (p8)  --  (p10);
           \draw (p3)  --  (p4);
           \draw (p5)  --  (p6);
           \draw (p7)  --  (p8);
           \draw (p9)  --  (p10);
           
           \draw (p1)  --  (p11);
        \draw (p2)  --  (p12);
      \draw (p11)  --  (p3);
        \draw (p7)  --  (p12);
     \draw (p11)  --  (p12);
           \draw (p15)  --  (p5);
           \draw (p15)  --  (p3);
           \draw (p15)  --  (p11);
           \draw (p16)  --  (p9);
           \draw (p16)  --  (p7);
           \draw (p16)  --  (p12);
           
            \draw (p1)  --  (p11r);
        \draw (p2)  --  (p12r);
      \draw (p11r)  --  (p4);
        \draw (p8)  --  (p12r);
     \draw (p11r)  --  (p12r);
           \draw (p15r)  --  (p6);
           \draw (p15r)  --  (p4);
           \draw (p15r)  --  (p11r);
           \draw (p16r)  --  (p10);
           \draw (p16r)  --  (p8);
           \draw (p16r)  --  (p12r);
           
          \draw (p16)  --  (p15);
           \draw (p16r)  --  (p15r);  
           
    \draw[dashed, thin](0,-1.9)--(0,1.9); 
 \draw[dashed, thin](-2.5,0)--(2.5,0); 
 \node [draw=white, fill=white] (b) at (0,-2.2) {(b)};
\end{tikzpicture}

\end{center}
\vspace{-0.6cm} \caption{Frameworks with $\mathcal{C}_{2v}$ symmetry discussed in Example~\ref{ex:c2v}. The framework in (a) has $k=-2$ and $e_{\sigma_h}=5$, $e_{\sigma_v}=3$.  It follows that it has $2$ fully-symmetric self-stresses and an anti-symmetric self-stress with respect to $\sigma_v$.
The framework in (b) has $k=-4$ and $e_{\sigma_h}=e_{\sigma_v}=5$. So this framework has $3$ fully-symmetric self-stresses and an anti-symmetric self-stress for each mirror. Note that both frameworks have an $A_2$-symmetric symmetry-detectable mechanism since $\alpha_2>0$.} \label{fig:2cvstresses}
\end{figure}

\item Suppose $k>0$. We again focus on the case when $k$ is even. (The other cases are analogous.) In this case we have $\alpha_2\geq 0$ and $\beta_2\geq 0$ for any values of $e_{\sigma_h}$ and $e_{\sigma_v}$. We also have $\beta_1\geq \alpha_1$, so if $\alpha_1\geq 0$, or equivalently,  $e_{\sigma_h}+e_{\sigma_v}-2\leq k$, then $\beta_1\geq 0$ and we only detect the $k$ mechanisms predicted by the standard Maxwell rule. So suppose $\alpha_1<0$. If $\beta_1\leq 0$, or equivalently, $e_{\sigma_h}-e_{\sigma_v}\geq k$, then we detect $-\alpha_1-\beta_1=\frac{-k+e_{\sigma_h}-1}{2}$ self-stresses -- the same amount as with a $\mathcal{C}_s$ analysis with the $\sigma_h$ mirror.
If $\alpha_1<0$ and $\beta_1> 0$, or equivalently, $e_{\sigma_h}+e_{\sigma_v}-2>k> e_{\sigma_h}-e_{\sigma_v}$, then  $\beta_2>0$ and $\alpha_2>0$, and we detect  $-\alpha_1$ fully-symmetric self-stresses, which is more than we detect with a $\mathcal{C}_s$ analysis.

\item For a fixed value of $k$ we increase the number of fully-symmetric self-stresses (and $A_2$-symmetric mechanisms) by increasing the total number of bars that are unshifted by a mirror, i.e., by increasing $e_{\sigma_h}+e_{\sigma_v}$. To increase the number of $B_1$-symmetric self-stresses (i.e., self-stresses that are anti-symmetric with respect to $\sigma_v$) we need to make $e_{\sigma_v}$ small in comparison to $e_{\sigma_h}$. This is consistent with what we observed for frameworks with $\mathcal{C}_s$ symmetry. The framework in Figure~\ref{fig:2cvstresses}(a) illustrates this.

As we observed in (i), by choosing $e_{\sigma_h}+e_{\sigma_v}$ sufficiently large and by keeping the difference between  $e_{\sigma_h}$ and $e_{\sigma_v}$ suitably small, we may obtain self-stresses of symmetry types $A_1, B_1$ and $B_2$ (and mechanisms of type $A_2$). See  Figure~\ref{fig:2cvstresses}(b) for an example. Note that such a distribution of self-stresses is particularly useful for the construction of gridshells.
\end{enumerate}

\begin{example} \label{ex:c2v}
Figure~\ref{fig:2cvstresses} shows two examples of  frameworks with $\mathcal{C}_{2v}$ symmetry. 
The framework in (a)  has $e=2v-1=31$, so $k=-2$, and $e_{\sigma_h}=5$, $e_{\sigma_v}=3$.
Thus, by Equation~(\ref{eq:c2vcase2}), we have $\Gamma(m)-\Gamma(s)= -2A_1+A_2-B_1$. So this framework has $2$ fully-symmetric self-stresses and an anti-symmetric self-stress with respect to $\sigma_v$. A $\mathcal{C}_s$ analysis with the reflection $\sigma_h$ also finds three self-stresses, all of which are fully-symmetric with respect to $\sigma_h$: $\Gamma(m)-\Gamma(s)=-3A'+A''$.

The framework in (b) has $e=2v+1=37$, so $k=-4$, and $e_{\sigma_h}=e_{\sigma_v}=5$. Thus, by Equation~(\ref{eq:c2vcase2}), we have $\Gamma(m)-\Gamma(s)= -3A_1+A_2-B_1-B_2$. So this framework has $3$ fully-symmetric self-stresses and an anti-symmetric self-stress for each mirror. Note that a $\mathcal{C}_s$ analysis (with either mirror) only detects 4 self-stresses: $\Gamma(m)-\Gamma(s)=-4A'$.
\end{example}

\subsection{Dihedral symmetry $\mathcal{C}_{nv}$, $n\geq 3$} \label{sec:dih}

In this section we consider dihedral symmetry groups of order at least 6. For simplicity, we will focus on the groups $\mathcal{C}_{3v}$ and $\mathcal{C}_{4v}$, but the groups $\mathcal{C}_{nv}$ with $n\geq 5$ can be analysed analogously. The characters of the irreducible representations of $\mathcal{C}_{3v}$  and  $\mathcal{C}_{4v}$ are shown in Table~\ref{tab:char}. Note that $\mathcal{C}_{3v}$  and  $\mathcal{C}_{4v}$ are of order 6 and 8, respectively. However, since for every element of the group that lies in the same conjugacy class we obtain the same trace, the tables only have one column for each conjugacy class of the group. The number of elements in each conjugacy class is indicated by the coefficient in front of the element that represents this conjugacy class in the character table. 

For example, $2C_3$ in the $\mathcal{C}_{3v}$ table stands for the rotations $C_3$  and $C_3^2$ about the origin by  $\frac{2\pi}{3}$ and $\frac{4\pi}{3}$, respectively, which lie in the same conjugacy class of $\mathcal{C}_{3v}$. For $\mathcal{C}_{4v}$, $2\sigma_v$ stands for the reflections in the vertical and horizontal mirrors, and $2\sigma_d$ stands for the reflections in the two diagonal mirrors.

\begin{table}[htp]
\begin{center}
\begin{tabular}{c||c|c|c}
$\mathcal{C}_{3v}$   &   $E$  &  $2C_{3}$   & $3\sigma$ \\\hline\hline
$A_{1}$  &    1  &  1 &  1  \\\hline
$A_{2}$  &    1  &  1 &  -1  \\\hline
$E$  &    2  &  -1 &  0  \\
\end{tabular}\quad\quad
\begin{tabular}{c||c|c|c|c|c}
$\mathcal{C}_{4v}$   &   $E$  &  $2C_{4}$   & $C_{2}$   &    $2\sigma_{v}$ &  $2\sigma_{d}$ \\\hline\hline
$A_{1}$  &    1  &  1 &  1  &  1 & 1   \\\hline
$A_{2}$  &    1  &  1 &  1  &  -1 &    -1   \\\hline
$B_{1}$  &    1  &  -1 &  1  &  1 &    -1  \\\hline
$B_{2}$  &    1  &  -1 &  1  &  -1 &    1  \\\hline
$E$  &    2  &  0 &  -2  &  0 & 0  \\
\end{tabular}
\caption{The irreducible characters of $\mathcal{C}_{3v}$ and $\mathcal{C}_{4v}$.} \label{tab:char}\end{center}
\end{table}

Using the same approch as above, we will  derive formulas for $\Gamma(m)-\Gamma(s)$ for the groups $\mathcal{C}_{3v}$ and $\mathcal{C}_{4v}$. We will also make some observations arising from these formulas in each case. However, since these analyses are similar to the one we have done for $\mathcal{C}_{2v}$, we will keep this discussion fairly succinct by focusing on the cases when $k\leq 0$ and $v_c=0$.

\subsubsection{The group $\mathcal{C}_{3v}$}

For a planar framework with $\mathcal{C}_{3v}$ symmetry, we have $v_c=0$ or $1$ and $e_3=0$. Suppose the framework has freedom number $k$. Then, by Table~\ref{tab:2D}, we have
\begin{equation*}\Gamma(m)-\Gamma(s)=(k,-v_c,-e_{\sigma}+1).\label{eq:c3v}\end{equation*}
Using Equation~(\ref{eq:coef}) we then obtain:
\begin{equation}\label{eq:c3vfinal}
\Gamma(m)-\Gamma(s)= \begin{cases}
\frac{k-3e_{\sigma}+3}{6}A_1+\frac{k+3e_{\sigma}-3}{6}A_2+ \frac{k}{3}E, &  \textrm{ if } v_c=0 \\
\frac{k-3e_{\sigma}+1}{6}A_1+\frac{k+3e_{\sigma}-5}{6}A_2+ \frac{k+1}{3}E, &  \textrm{ if } v_c=1\\
\end{cases}
\end{equation}
Note that if $k$ is even, then $e$ is odd and hence $e_{\sigma}$ is odd. Similarly, if $k$ is odd, then $e$ is even and hence $e_{\sigma}$ is even. Also, $k$ is divisible by $3$ if and only if $v_c=0$, and $k+1$ is divisible by $3$ if and only if $v_c=1$. 

Some observations arising from Equation~(\ref{eq:c3vfinal}) are:
\begin{enumerate}
\item We focus on the case $v_c=0$, as the case $v_c=1$ is analogous. Suppose $k\leq 0$ and $k$ is even. Note that this implies that if $k$ is non-zero, we have $k\leq -6$. We will denote the coefficients of $A_i$ by $\alpha_i$ for $i=1,2$, and the coefficient of $E$ by $\epsilon$. We have $\alpha_1+\alpha_2+2\epsilon=k$ (since $E$ is the character of a $2$-dimensional representation), and $\alpha_1\leq 0$ and $\epsilon\leq 0$ for any value of $e_{\sigma}$. 
If $\alpha_2\leq 0$, or equivalently $e_\sigma\leq \frac{-k+3}{3}$,    then we only find the $-k$ self-stresses predicted by the standard Maxwell rule.  If $\alpha_2>0$, that is, $e_\sigma>\frac{-k+3}{3}$, then we detect  $-\alpha_1-2\epsilon=\frac{-5k+3e_\sigma-3}{6}>-k$ self-stresses, which is more than we detect with a $\mathcal{C}_s$ analysis, provided that $k\neq 0$ (recall Section~\ref{sec:mirror}). 
We may draw similar conclusions if $k\leq 0$ and $k$ is odd.

\item In the special case of $k=0$, we must have $v_c=0$, and there are no symmetry-detectable self-stresses or mechanisms if $e_\sigma=1$. In this case the framework is \emph{conjectured} to be isostatic for any `generic' positions of the vertices \cite{cfgsw}.  If $e_\sigma\geq 3$, then we find $\frac{e_\sigma-1}{2}$ symmetry-detectable fully-symmetric self-stresses. 

\item Analogous to the $\mathcal{C}_s$ situation, increasing $e_{\sigma}$ while keeping $k$ fixed increases the number of fully-symmetric self-stresses.  The number of $E$-symmetric self-stresses only depends on $k$. 
\end{enumerate}

\subsubsection{The group $\mathcal{C}_{4v}$}

For a planar framework with $\mathcal{C}_{4v}$ symmetry, we have $v_c=0$ or $1$ and $e_4=e_2=0$. Suppose the framework has freedom number $k$. Then, by Table~\ref{tab:2D}, we have
\begin{equation*}\Gamma(m)-\Gamma(s)=(k,-1, -2v_c+1 ,-e_{\sigma_v}+1, -e_{\sigma_d}+1).\label{eq:c3v}\end{equation*}
Using Equation~(\ref{eq:coef}) we then obtain the following expressions for $\Gamma(m)-\Gamma(s)$:

\begin{itemize}
\item For $v_c=0$ we obtain:
{
\medmuskip=0mu
\thinmuskip=0mu
\thickmuskip=0mu
\nulldelimiterspace=0pt
\scriptspace=-1pt
\begin{equation*}\label{eq:c4vfinal0}
\hspace{-0.7cm}\resizebox{0.98\hsize}{!}{$\Gamma(m)-\Gamma(s)=\frac{k-2e_{\sigma_v}-2e_{\sigma_d}+3}{8}A_1+\frac{k+2e_{\sigma_v}+2e_{\sigma_d}-5}{8}A_2+ \frac{k-2e_{\sigma_v}+2e_{\sigma_d}+3}{8}B_1+\frac{k+2e_{\sigma_v}-2e_{\sigma_d}+3}{8}B_2+ \frac{k-1}{4}E$}
\end{equation*}
}
\item For $v_c=1$ we obtain:
{
\medmuskip=0mu
\thinmuskip=0mu
\thickmuskip=0mu
\nulldelimiterspace=0pt
\scriptspace=-1pt
\begin{equation*}\label{eq:c4vfinal1}
\hspace{-0.7cm}\resizebox{0.98\hsize}{!}{$\Gamma(m)-\Gamma(s)=\frac{k-2e_{\sigma_v}-2e_{\sigma_d}+1}{8}A_1+\frac{k+2e_{\sigma_v}+2e_{\sigma_d}-7}{8}A_2
+ \frac{k-2e_{\sigma_v}+2e_{\sigma_d}+1}{8}B_1+\frac{k+2e_{\sigma_v}-2e_{\sigma_d}+1}{8}B_2
+ \frac{k+1}{4}E$}
\end{equation*}
}
\end{itemize}

Note that   $e_4=e_2=0$ implies that  $e_{\sigma_v}$ and $e_{\sigma_d}$ are even. Hence $e$ is even  and  $k$ is odd.  Also, $k-1$ is divisible by $4$ if and only if $v_c=0$, and $k+1$ is divisible by $4$ if and only if $v_c=1$. 

\begin{figure}[htp]
\begin{center}
\begin{tikzpicture}[very thick,scale=0.7]
\tikzstyle{every node}=[circle, draw=black, fill=white, inner sep=0pt, minimum width=4pt];
   
       \path (-0.4,0.4) node (p1)  {} ;
    \path (-0.4,1.2) node (p2)  {} ;
            \path (-0.4,2) node (p3)  {} ;
     \path (0.4,0.4) node (p1r)  {} ;
    \path (0.4,1.2) node (p2r)  {} ;
            \path (0.4,2) node (p3r)  {} ;
           \path (-0.4,-0.4) node (p4)  {} ;
    \path (-0.4,-1.2) node (p5)  {} ;
            \path (-0.4,-2) node (p6)  {} ;
     \path (0.4,-0.4) node (p4r)  {} ;
    \path (0.4,-1.2) node (p5r)  {} ;
            \path (0.4,-2) node (p6r)  {} ;
 \path (-1.2,1.2) node (p7)  {} ;
    \path (-1.2,0.4) node (p8)  {} ;
            \path (-1.2,-0.4) node (p9)  {} ;
 \path (-1.2,-1.2) node (p10)  {} ;
 \path (1.2,1.2) node (p7r)  {} ;
    \path (1.2,0.4) node (p8r)  {} ;
            \path (1.2,-0.4) node (p9r)  {} ;
 \path (1.2,-1.2) node (p10r)  {} ;

      \path (-1.8,1.8) node (p11)  {} ;
    \path (-2,0.4) node (p12)  {} ;
            \path (-2,-0.4) node (p13)  {} ;
 \path (-1.8,-1.8) node (p14)  {} ;
 \path (1.8,1.8) node (p11r)  {} ;
    \path (2,0.4) node (p12r)  {} ;
            \path (2,-0.4) node (p13r)  {} ;
 \path (1.8,-1.8) node (p14r)  {} ;
 
  \draw (p1)  --  (p2);
        \draw (p2)  --  (p3);
    \draw (p1)  --  (p4);
        \draw (p4)  --  (p5);
            \draw (p5)  --  (p6);
        \draw (p7)  --  (p8);
     \draw (p8)  --  (p9);
        \draw (p9)  --  (p10);
       \draw (p11)  --  (p12);
        \draw (p12)  --  (p13);
     \draw (p13)  --  (p14);
        \draw (p3)  --  (p11);  
         \draw (p3)  --  (p7);
        \draw (p7)  --  (p11);
     \draw (p2)  --  (p7);
        \draw (p1)  --  (p8);
         \draw (p4)  --  (p9);
        \draw (p5)  --  (p10);
     \draw (p6)  --  (p10);
        \draw (p6)  --  (p14);
         \draw (p14)  --  (p10);
        \draw (p10)  --  (p13);
     \draw (p13)  --  (p9);
        \draw (p8)  --  (p12);
        \draw (p7)  --  (p12);   
           
     \draw (p1r)  --  (p2r);
        \draw (p2r)  --  (p3r);
    \draw (p1r)  --  (p4r);
        \draw (p4r)  --  (p5r);
            \draw (p5r)  --  (p6r);
        \draw (p7r)  --  (p8r);
     \draw (p8r)  --  (p9r);
        \draw (p9r)  --  (p10r);
       \draw (p11r)  --  (p12r);
        \draw (p12r)  --  (p13r);
     \draw (p13r)  --  (p14r);
        \draw (p3r)  --  (p11r);  
         \draw (p3r)  --  (p7r);
        \draw (p7r)  --  (p11r);
     \draw (p2r)  --  (p7r);
        \draw (p1r)  --  (p8r);
         \draw (p4r)  --  (p9r);
        \draw (p5r)  --  (p10r);
     \draw (p6r)  --  (p10r);
        \draw (p6r)  --  (p14r);
         \draw (p14r)  --  (p10r);
        \draw (p10r)  --  (p13r);
     \draw (p13r)  --  (p9r);
        \draw (p8r)  --  (p12r);
        \draw (p7r)  --  (p12r);        
           
         \draw (p3)  --  (p3r);
        \draw (p2)  --  (p2r);   
           \draw (p1)  --  (p1r);
        \draw (p4)  --  (p4r); 
           \draw (p5)  --  (p5r);
        \draw (p6)  --  (p6r);

    \draw[dashed, thin](0,-2.5)--(0,2.5); 
 \draw[dashed, thin](-2.5,0)--(2.5,0);
  \draw[dashed, thin](-2.5,-2.5)--(2.5,2.5); 
 \draw[dashed, thin](2.5,-2.5)--(-2.5,2.5); 
  
 \node [draw=white, fill=white] (b) at (0,-3.3) {(a)};
\end{tikzpicture}
\hspace{1cm}
\begin{tikzpicture}[very thick,scale=0.7]
\tikzstyle{every node}=[circle, draw=black, thick, fill=white, inner sep=0pt, minimum width=3pt];
             
         \path (-0.4,0.8) node (p1)  {} ;
    \path (-0.4,1.4) node (p2)  {} ;
            \path (-0.4,2) node (p3)  {} ;
             \path (0,2.6) node (p4)  {} ;
     \path (0.4,0.8) node (p11)  {} ;
    \path (0.4,1.4) node (p22)  {} ;
            \path (0.4,2) node (p33)  {} ;
           \path (-0.4,-0.8) node (p1u)  {} ;
    \path (-0.4,-1.4) node (p2u)  {} ;
            \path (-0.4,-2) node (p3u)  {} ;
     \path (0.4,-0.8) node (p11u)  {} ;
    \path (0.4,-1.4) node (p22u)  {} ;
            \path (0.4,-2) node (p33u)  {} ;     
              \path (0,-2.6) node (p4u)  {} ;
              
          \path (-0.8,0.4) node (p1l)  {} ;
    \path (-1.4,0.4) node (p2l)  {} ;
            \path (-2,0.4) node (p3l)  {} ;
             \path (-2.6,0) node (p4l)  {} ;
    \path (-0.8,-0.4) node (p11l)  {} ;
    \path (-1.4,-0.4) node (p22l)  {} ;
            \path (-2,-0.4) node (p33l)  {} ;
          \path (0.8,0.4) node (p1r)  {} ;
    \path (1.4,0.4) node (p2r)  {} ;
            \path (2,0.4) node (p3r)  {} ;
             \path (2.6,0) node (p4r)  {} ;
    \path (0.8,-0.4) node (p11r)  {} ;
    \path (1.4,-0.4) node (p22r)  {} ;
            \path (2,-0.4) node (p33r)  {} ;
           
          \path (-0.9,1.9) node (q1)  {} ;
    \path (-2.1,2.1) node (m1)  {} ;
            \path (-1.9,0.9) node (q1l)  {} ;        
               \path (-0.9,2.5) node (b1)  {} ;
                \path (-2.5,0.9) node (b1l)  {} ;
               
                  \path (0.9,1.9) node (q1r)  {} ;
    \path (2.1,2.1) node (m1r)  {} ;
            \path (1.9,0.9) node (q1lr)  {} ;        
               \path (0.9,2.5) node (b1r)  {} ;
               \path (2.5,0.9) node (b1lr)  {} ;
               
              \path (-0.9,-1.9) node (q1u)  {} ;
    \path (-2.1,-2.1) node (m1u)  {} ;
            \path (-1.9,-0.9) node (q1lu)  {} ;        
               \path (-0.9,-2.5) node (b1u)  {} ;
                \path (-2.5,-0.9) node (b1luu)  {} ;
               
                  \path (0.9,-1.9) node (q1ru)  {} ;
    \path (2.1,-2.1) node (m1ru)  {} ;
            \path (1.9,-0.9) node (q1lru)  {} ;        
               \path (0.9,-2.5) node (b1ru)  {} ;
               \path (2.5,-0.9) node (b1lu)  {} ;

      \draw (p1)  --  (p2);
        \draw (p2)  --  (p3);
    \draw (p3)  --  (p4);
        \draw (p11)  --  (p22);
            \draw (p22)  --  (p33);
     \draw (p33)  --  (p4);
     
     \draw (p1u)  --  (p2u);
        \draw (p2u)  --  (p3u);
    \draw (p3u)  --  (p4u);
        \draw (p11u)  --  (p22u);
            \draw (p22u)  --  (p33u);
     \draw (p33u)  --  (p4u);
     
     \draw (p1l)  --  (p2l);
        \draw (p2l)  --  (p3l);
    \draw (p3l)  --  (p4l);
        \draw (p11l)  --  (p22l);
            \draw (p22l)  --  (p33l);
     \draw (p33l)  --  (p4l);
     
     \draw (p1r)  --  (p2r);
        \draw (p2r)  --  (p3r);
    \draw (p3r)  --  (p4r);
        \draw (p11r)  --  (p22r);
            \draw (p22r)  --  (p33r);
     \draw (p33r)  --  (p4r);
     
        \draw (p1)  --  (p11);
        \draw (p11)  --  (p1r);
    \draw (p1r)  --  (p11r);
        \draw (p11r)  --  (p11u);
            \draw (p11u)  --  (p1u);
     \draw (p1u)  --  (p11l);  
             \draw (p11l)  --  (p1l);
     \draw (p1l)  --  (p1);  
          
          \draw (p2)  --  (p22);
        \draw (p22)  --  (p2r);
    \draw (p2r)  --  (p22r);
        \draw (p22r)  --  (p22u);
            \draw (p22u)  --  (p2u);
     \draw (p2u)  --  (p22l);  
             \draw (p22l)  --  (p2l);
     \draw (p2l)  --  (p2);
     
    \draw (p3)  --  (p33);
        \draw (p33)  --  (q1r);
    \draw (q1r)  --  (q1lr);
        \draw (q1lr)  --  (p3r);
            \draw (p3r)  --  (p33r);
     \draw (p33r)  --  (q1lru);  
             \draw (q1lru)  --  (q1ru);
        \draw (q1ru)  --  (p33u);
    \draw (p33u)  --  (p3u);
        \draw (p3u)  --  (q1u);
            \draw (q1u)  --  (q1lu);
     \draw (q1lu)  --  (p33l);  
      \draw (p33l)  --  (p3l);
        \draw (p3l)  --  (q1l);
            \draw (q1l)  --  (q1);
     \draw (q1)  --  (p3); 
     
      \draw (p22)  --  (q1r);
        \draw (p2)  --  (q1);
      \draw (p2r)  --  (q1lr);
        \draw (p22r)  --  (q1lru);
        
         \draw (p22u)  --  (q1ru);
        \draw (p2u)  --  (q1u);
      \draw (p22l)  --  (q1lu);
        \draw (p2l)  --  (q1l);
        
         \draw (p4)  --  (b1r);
        \draw (b1r)  --  (m1r);
      \draw (m1r)  --  (b1lr);
        \draw (b1lr)  --  (p4r);
                 \draw (p4r)  --  (b1lu);
        \draw (b1lu)  --  (m1ru);
      \draw (m1ru)  --  (b1ru);
        \draw (b1ru)  --  (p4u);
          \draw (p4u)  --  (b1u);
        \draw (b1u)  --  (m1u);
      \draw (m1u)  --  (b1luu);
        \draw (b1luu)  --  (p4l);
            \draw (p4l)  --  (b1l);
        \draw (b1l)  --  (m1);
      \draw (m1)  --  (b1);
        \draw (b1)  --  (p4); 
             
         \draw (m1r)  --  (q1r);
        \draw (m1r)  --  (q1lr);     
             \draw (m1ru)  --  (q1lru);
        \draw (m1ru)  --  (q1ru); 
        \draw (m1u)  --  (q1u);
        \draw (m1u)  --  (q1lu); 
        \draw (m1)  --  (q1);
        \draw (m1)  --  (q1l); 
        
         \draw (b1r)  --  (p33);
        \draw (b1r)  --  (q1r);
           \draw (b1lr)  --  (q1lr);
        \draw (b1lr)  --  (p3r);
                \draw (b1lu)  --  (p33r);
        \draw (b1lu)  --  (q1lru);
                \draw (b1ru)  --  (p33u);
        \draw (b1ru)  --  (q1ru);
                \draw (b1u)  --  (p3u);
        \draw (b1u)  --  (q1u);
                \draw (b1luu)  --  (p33l);
        \draw (b1luu)  --  (q1lu);
           \draw (b1l)  --  (p3l);
        \draw (b1l)  --  (q1l);
           \draw (b1)  --  (p3);
        \draw (b1)  --  (q1);
             
  \draw[dashed, thin](0,-3)--(0,3); 
 \draw[dashed, thin](-3.2,0)--(3.2,0);
  \draw[dashed, thin](-2.5,-2.5)--(2.5,2.5); 
 \draw[dashed, thin](2.5,-2.5)--(-2.5,2.5);
 \node [draw=white, fill=white] (b) at (0,-3.3) {(b)};
\end{tikzpicture}

\end{center}
\vspace{-0.6cm} \caption{Frameworks with $\mathcal{C}_{4v}$ symmetry discussed in Example~\ref{ex:c4v}. The framework in (a) has $k=-3$ and $e_{\sigma_v}=6$, $e_{\sigma_d}=2$.  It follows that it has two fully-symmetric self-stresses, a $B_1$-symmetric self-stress and two $E$-symmetric self-stresses, as well as an $A_2$-symmetric mechanism.
The framework in (b) has $k=-11$ and $e_{\sigma_v}=e_{\sigma_d}=6$.  A $\mathcal{C}_{4v}$ analysis finds 12 self-stresses, whereas a $\mathcal{C}_{2v}$ analysis only finds the 11 self-stresses predicted by the $k=-11$ count.} \label{fig:4cvstresses}
\end{figure}
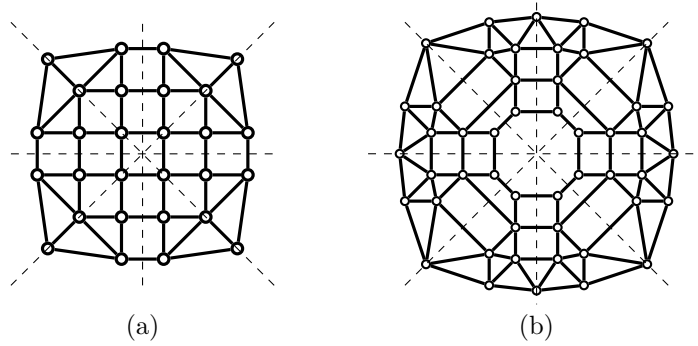

Some observations arising from these expressions for $\Gamma(m)-\Gamma(s)$ are: 
\begin{enumerate}
\item We focus on the case $v_c=0$, as the case $v_c=1$ is analogous. Suppose $k< 0$ and  $e_{\sigma_v}\geq e_{\sigma_d}$. We again denote the coefficients of $A_i$ and $B_i$ by $\alpha_i$ and $\beta_i$, respectively, for $i=1,2$, and the coefficient of $E$ by $\epsilon$. We have $\alpha_1+\alpha_2+\beta_1+\beta_2+2\epsilon=k$ (since $E$ is the character of a $2$-dimensional representation), and $\alpha_1,\beta_1,\epsilon\leq 0$ for any values of $e_{\sigma_v}$ and $e_{\sigma_d}$.

 Suppose that $e_{\sigma_d}\geq 2$. (The case when $e_{\sigma_d}=0$ is similar but less relevant for practical applications, since it forces the form diagram to be quite special.) We have $\alpha_2\geq\beta_2$. So if $\alpha_2\leq 0$,  then $\beta_2\leq 0$, and we only find the $-k$ self-stresses predicted by the standard Maxwell rule. 
 
 So suppose $\alpha_2> 0$ or equivalently $e_{\sigma_v}+e_{\sigma_d}>\frac{-k+5}{2}$. Then, if $\beta_2\geq 0$ or equivalently $e_{\sigma_v}-e_{\sigma_d}\geq\frac{-k-3}{2}$, we detect $-\alpha_1-\beta_1-2\epsilon$ self-stresses. This is  the same amount of self-stresses as we detect with a $\mathcal{C}_{2v}$ analysis (as is easily verified by considering Equation~(\ref{eq:c2vcase1}) in Section~\ref{sec:c2v}), but it is more than we detect with a $\mathcal{C}_s$ analysis (recall Section~\ref{sec:mirror}). See Figure~\ref{fig:4cvstresses}(a) for an example. 
 
 If $\alpha_2>0$ and we also have $\beta_2< 0$ or equivalently $e_{\sigma_v}-e_{\sigma_d}<\frac{-k-3}{2}$,  then we detect $-\alpha_1-\beta_1-\beta_2-2\epsilon$ self-stresses. In this case we find more self-stresses than with a $\mathcal{C}_{2v}$ analysis. 
 See Figure~\ref{fig:4cvstresses}(b) for an example.

\item If we fix $k$, then, analogously to the $\mathcal{C}_{2v}$ situation, we increase the number of fully-symmetric self-stresses (and $A_2$-symmetric mechanisms) by increasing the total number of bars that are unshifted by a mirror, i.e., by increasing $e_{\sigma_v}+e_{\sigma_d}$. To increase the number of $B_1$-symmetric self-stresses (i.e., self-stresses that are anti-symmetric with respect to $\sigma_d$) we need to 
 make $e_{\sigma_d}$ small in comparison to $e_{\sigma_v}$.
As observed above, by choosing  $e_{\sigma_v}+e_{\sigma_d}$  sufficiently large  and by keeping the difference between  $e_{\sigma_v}$ and $e_{\sigma_d}$ suitably small, we may obtain self-stresses of symmetry types $A_1, B_1$, $B_2$ and $E$ (and mechanisms of type $A_2$). Finally, note  that the number of $E$-symmetric self-stresses only depends on $k$.
\end{enumerate}

\begin{example} \label{ex:c4v}  Figure~\ref{fig:4cvstresses} shows two examples of  frameworks with $\mathcal{C}_{4v}$ symmetry. 
The framework in (a)  has $e=2v=56$, so $k=-3$. We also have $v_c=0$  and $e_{\sigma_v}=6$, $e_{\sigma_d}=2$.
Thus, we have $\Gamma(m)-\Gamma(s)= -2A_1+A_2-B_1+B_2-E$. So this framework has at least $5$ self-stresses, including $2$ fully-symmetric self-stresses and an anti-symmetric self-stress with respect to $\sigma_d$. A $\mathcal{C}_{2v}$ analysis with the vertical and horizontal mirror also finds 5 self-stresses: $\Gamma(m)-\Gamma(s)=-3A_1+2A_2-B_1-B_2$. However, a $\mathcal{C}_s$ analysis (with $\sigma_v$) only finds 4.

The framework in (b) has $e=2v+8=104$, so $k=-11$. We also have $v_c=0$ and $e_{\sigma_v}=e_{\sigma_d}=6$. Thus, we have $\Gamma(m)-\Gamma(s)= -4A_1+A_2-B_1-B_2-3E$. So  this framework has at least 12 self-stresses, including $4$ fully-symmetric self-stresses and an anti-symmetric self-stress for each pair of perpendicular mirrors. It also has a symmetry-detectable $A_2$-symmetric mechanism.
Note that a $\mathcal{C}_{2v}$ analysis of this framework (with either pair of perpendicular mirrors) only detects 11 self-stresses: $\Gamma(m)-\Gamma(s)=-5A_1-3B_1-3B_2$. A similar example for the case when $v_c=1$ is shown in Example~\ref{ex:spider}.
\end{example}

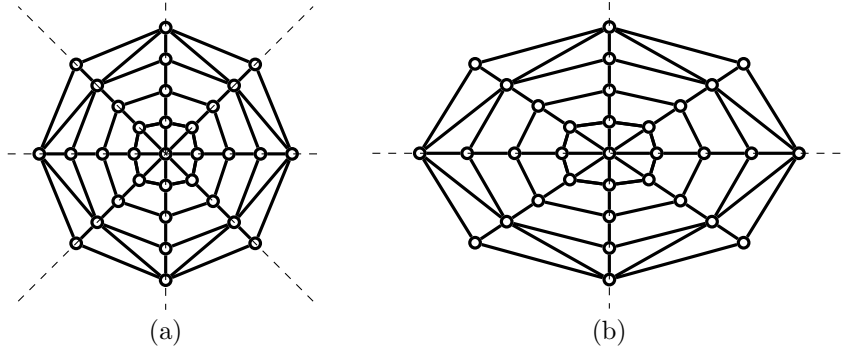
\begin{figure}[htp]
\begin{center}
\begin{tikzpicture}[very thick,scale=0.7]
\tikzstyle{every node}=[circle, draw=black, fill=white, inner sep=0pt, minimum width=4pt];
   
       \path (0,0) node (p0)  {} ;
       
    \path (0,0.6) node (p1)  {} ;
            \path (0,1.2) node (p2)  {} ;
     \path (0,1.8) node (p3)  {} ;
    \path (0,2.4) node (p4)  {} ;
            
         \path (0,-0.6) node (p1u)  {} ;
            \path (0,-1.2) node (p2u)  {} ;
     \path (0,-1.8) node (p3u)  {} ;
    \path (0,-2.4) node (p4u)  {} ;    
            
             \path (0.6,0) node (p1r)  {} ;
            \path (1.2,0) node (p2r)  {} ;
     \path (1.8,0) node (p3r)  {} ;
    \path (2.4,0) node (p4r)  {} ;
    
       \path (-0.6,0) node (p1l)  {} ;
            \path (-1.2,0) node (p2l)  {} ;
     \path (-1.8,0) node (p3l)  {} ;
    \path (-2.4,0) node (p4l)  {} ;

     \path (0.5,0.5) node (d1)  {} ;
            \path (0.9,0.9) node (d2)  {} ;
     \path (1.3,1.3) node (d3)  {} ;
    \path (1.7,1.7) node (d4)  {} ;
     
        \path (-0.5,0.5) node (d1l)  {} ;
            \path (-0.9,0.9) node (d2l)  {} ;
     \path (-1.3,1.3) node (d3l)  {} ;
    \path (-1.7,1.7) node (d4l)  {} ;
     
      \path (-0.5,-0.5) node (d1lu)  {} ;
            \path (-0.9,-0.9) node (d2lu)  {} ;
     \path (-1.3,-1.3) node (d3lu)  {} ;
    \path (-1.7,-1.7) node (d4lu)  {} ;
     
      \path (0.5,-0.5) node (d1r)  {} ;
            \path (0.9,-0.9) node (d2r)  {} ;
     \path (1.3,-1.3) node (d3r)  {} ;
    \path (1.7,-1.7) node (d4r)  {} ;
     
  \draw (p0)  --  (p1);
      \draw (p1)  --  (p2);
      \draw (p2)  --  (p3);
          \draw (p3)  --  (p4);
    
    \draw (p0)  --  (p1r);
      \draw (p1r)  --  (p2r);
      \draw (p2r)  --  (p3r);
          \draw (p3r)  --  (p4r);
    
    \draw (p0)  --  (p1l);
      \draw (p1l)  --  (p2l);
      \draw (p2l)  --  (p3l);
          \draw (p3l)  --  (p4l);
    
    \draw (p0)  --  (p1u);
      \draw (p1u)  --  (p2u);
      \draw (p2u)  --  (p3u);
          \draw (p3u)  --  (p4u);
    
    \draw (p0)  --  (d1);
      \draw (d1)  --  (d2);
      \draw (d2)  --  (d3);
          \draw (d3)  --  (d4);
    
      \draw (p0)  --  (d1l);
      \draw (d1l)  --  (d2l);
      \draw (d2l)  --  (d3l);
          \draw (d3l)  --  (d4l);
           
          \draw (p0)  --  (d1lu);
      \draw (d1lu)  --  (d2lu);
      \draw (d2lu)  --  (d3lu);
          \draw (d3lu)  --  (d4lu); 
           
           \draw (p0)  --  (d1r);
      \draw (d1r)  --  (d2r);
      \draw (d2r)  --  (d3r);
          \draw (d3r)  --  (d4r);
          
        \draw (p1)  --  (d1l);
      \draw (d1l)  --  (p1l);
      \draw (p1l)  --  (d1lu);
          \draw (d1lu)  --  (p1u);  
          \draw (p1u)  --  (d1r);
      \draw (d1r)  --  (p1r);
      \draw (p1r)  --  (d1);
          \draw (d1)  --  (p1); 
          
             \draw (p1)  --  (d1l);
      \draw (d1l)  --  (p1l);
      \draw (p1l)  --  (d1lu);
          \draw (d1lu)  --  (p1u);  
          \draw (p1u)  --  (d1r);
      \draw (d1r)  --  (p1r);
      \draw (p1r)  --  (d1);
          \draw (d1)  --  (p1);

            \draw (p2)  --  (d2l);
      \draw (d2l)  --  (p2l);
      \draw (p2l)  --  (d2lu);
          \draw (d2lu)  --  (p2u);  
          \draw (p2u)  --  (d2r);
      \draw (d2r)  --  (p2r);
      \draw (p2r)  --  (d2);
          \draw (d2)  --  (p2); 
          
                \draw (p3)  --  (d3l);
      \draw (d3l)  --  (p3l);
      \draw (p3l)  --  (d3lu);
          \draw (d3lu)  --  (p3u);  
          \draw (p3u)  --  (d3r);
      \draw (d3r)  --  (p3r);
      \draw (p3r)  --  (d3);
          \draw (d3)  --  (p3); 
          
                \draw (p4)  --  (d4l);
      \draw (d4l)  --  (p4l);
      \draw (p4l)  --  (d4lu);
          \draw (d4lu)  --  (p4u);  
          \draw (p4u)  --  (d4r);
      \draw (d4r)  --  (p4r);
      \draw (p4r)  --  (d4);
          \draw (d4)  --  (p4); 
          
                \draw (p4)  --  (d3);
      \draw (p4)  --  (d3l);
      \draw (d3l)  --  (p4l);
          \draw (p4l)  --  (d3lu);  
          \draw (d3lu)  --  (p4u);
      \draw (p4u)  --  (d3r);
      \draw (d3r)  --  (p4r);
          \draw (p4r)  --  (d3); 
          
    \draw[dashed, thin](0,-3)--(0,3); 
 \draw[dashed, thin](-3,0)--(3,0);
  \draw[dashed, thin](-2.8,-2.8)--(2.8,2.8); 
 \draw[dashed, thin](2.8,-2.8)--(-2.8,2.8); 
  
 \node [draw=white, fill=white] (b) at (0,-3.4) {(a)};
\end{tikzpicture}
\hspace{0.4cm}
\begin{tikzpicture}[very thick,scale=0.7]
\tikzstyle{every node}=[circle, draw=black, fill=white, inner sep=0pt, minimum width=4pt];
               \path (0,0) node (p0)  {} ;
       
    \path (0,0.6) node (p1)  {} ;
            \path (0,1.2) node (p2)  {} ;
     \path (0,1.8) node (p3)  {} ;
    \path (0,2.4) node (p4)  {} ;
            
         \path (0,-0.6) node (p1u)  {} ;
            \path (0,-1.2) node (p2u)  {} ;
     \path (0,-1.8) node (p3u)  {} ;
    \path (0,-2.4) node (p4u)  {} ;    
            
             \path (0.9,0) node (p1r)  {} ;
            \path (1.8,0) node (p2r)  {} ;
     \path (2.7,0) node (p3r)  {} ;
    \path (3.6,0) node (p4r)  {} ;
    
       \path (-0.9,0) node (p1l)  {} ;
            \path (-1.8,0) node (p2l)  {} ;
     \path (-2.7,0) node (p3l)  {} ;
    \path (-3.6,0) node (p4l)  {} ;

     \path (0.75,0.5) node (d1)  {} ;
            \path (1.35,0.9) node (d2)  {} ;
     \path (1.95,1.3) node (d3)  {} ;
    \path (2.55,1.7) node (d4)  {} ;
     
        \path (-0.75,0.5) node (d1l)  {} ;
            \path (-1.35,0.9) node (d2l)  {} ;
     \path (-1.95,1.3) node (d3l)  {} ;
    \path (-2.55,1.7) node (d4l)  {} ;
     
      \path (-0.75,-0.5) node (d1lu)  {} ;
            \path (-1.35,-0.9) node (d2lu)  {} ;
     \path (-1.95,-1.3) node (d3lu)  {} ;
    \path (-2.55,-1.7) node (d4lu)  {} ;
     
      \path (0.75,-0.5) node (d1r)  {} ;
            \path (1.35,-0.9) node (d2r)  {} ;
     \path (1.95,-1.3) node (d3r)  {} ;
    \path (2.55,-1.7) node (d4r)  {} ;
     
  \draw (p0)  --  (p1);
      \draw (p1)  --  (p2);
      \draw (p2)  --  (p3);
          \draw (p3)  --  (p4);
    
    \draw (p0)  --  (p1r);
      \draw (p1r)  --  (p2r);
      \draw (p2r)  --  (p3r);
          \draw (p3r)  --  (p4r);
    
    \draw (p0)  --  (p1l);
      \draw (p1l)  --  (p2l);
      \draw (p2l)  --  (p3l);
          \draw (p3l)  --  (p4l);
    
    \draw (p0)  --  (p1u);
      \draw (p1u)  --  (p2u);
      \draw (p2u)  --  (p3u);
          \draw (p3u)  --  (p4u);
    
    \draw (p0)  --  (d1);
      \draw (d1)  --  (d2);
      \draw (d2)  --  (d3);
          \draw (d3)  --  (d4);
    
      \draw (p0)  --  (d1l);
      \draw (d1l)  --  (d2l);
      \draw (d2l)  --  (d3l);
          \draw (d3l)  --  (d4l);
           
          \draw (p0)  --  (d1lu);
      \draw (d1lu)  --  (d2lu);
      \draw (d2lu)  --  (d3lu);
          \draw (d3lu)  --  (d4lu); 
           
           \draw (p0)  --  (d1r);
      \draw (d1r)  --  (d2r);
      \draw (d2r)  --  (d3r);
          \draw (d3r)  --  (d4r);
          
        \draw (p1)  --  (d1l);
      \draw (d1l)  --  (p1l);
      \draw (p1l)  --  (d1lu);
          \draw (d1lu)  --  (p1u);  
          \draw (p1u)  --  (d1r);
      \draw (d1r)  --  (p1r);
      \draw (p1r)  --  (d1);
          \draw (d1)  --  (p1); 
          
             \draw (p1)  --  (d1l);
      \draw (d1l)  --  (p1l);
      \draw (p1l)  --  (d1lu);
          \draw (d1lu)  --  (p1u);  
          \draw (p1u)  --  (d1r);
      \draw (d1r)  --  (p1r);
      \draw (p1r)  --  (d1);
          \draw (d1)  --  (p1);

            \draw (p2)  --  (d2l);
      \draw (d2l)  --  (p2l);
      \draw (p2l)  --  (d2lu);
          \draw (d2lu)  --  (p2u);  
          \draw (p2u)  --  (d2r);
      \draw (d2r)  --  (p2r);
      \draw (p2r)  --  (d2);
          \draw (d2)  --  (p2); 
          
                \draw (p3)  --  (d3l);
      \draw (d3l)  --  (p3l);
      \draw (p3l)  --  (d3lu);
          \draw (d3lu)  --  (p3u);  
          \draw (p3u)  --  (d3r);
      \draw (d3r)  --  (p3r);
      \draw (p3r)  --  (d3);
          \draw (d3)  --  (p3); 
          
                \draw (p4)  --  (d4l);
      \draw (d4l)  --  (p4l);
      \draw (p4l)  --  (d4lu);
          \draw (d4lu)  --  (p4u);  
          \draw (p4u)  --  (d4r);
      \draw (d4r)  --  (p4r);
      \draw (p4r)  --  (d4);
          \draw (d4)  --  (p4); 
          
                \draw (p4)  --  (d3);
      \draw (p4)  --  (d3l);
      \draw (d3l)  --  (p4l);
          \draw (p4l)  --  (d3lu);  
          \draw (d3lu)  --  (p4u);
      \draw (p4u)  --  (d3r);
      \draw (d3r)  --  (p4r);
          \draw (p4r)  --  (d3);    
             
    \draw[dashed, thin](0,-3)--(0,3); 
 \draw[dashed, thin](-4.5,0)--(4.5,0); 
 \node [draw=white, fill=white] (b) at (0,-3.4) {(b)};
\end{tikzpicture}
\end{center}
\vspace{-0.6cm} \caption{Frameworks with $\mathcal{C}_{4v}$ symmetry. The framework in (a) has $k=-9$ and $e_{\sigma_v}=e_{\sigma_d}=8$.  A  $\mathcal{C}_{4v}$ analysis finds 11 self-stresses, as detailed in  Example~\ref{ex:spider}, whereas a $\mathcal{C}_{2v}$ analysis only finds 10.
The framework in (b) is obtained from the one in (a) by a horizontal stretch so that it only has $\mathcal{C}_{2v}$ symmetry.} \label{fig:spiderweb}
\end{figure}

\begin{example}  \label{ex:spider} Figure~\ref{fig:spiderweb}(a) shows another example of a planar framework with $\mathcal{C}_{4v}$ symmetry. This framework has $v_c=1$. For such frameworks, a similar analysis as in (i) shows that if $k<0$, $e_{\sigma_v}+e_{\sigma_d}>\frac{-k+7}{2}$ and $e_{\sigma_v}-e_{\sigma_d}<\frac{-k-1}{2}$, then we detect more self-stresses with a $\mathcal{C}_{4v}$ analysis than with a $\mathcal{C}_{2v}$ analysis. In particular, the framework is guaranteed to have fully-symmetric self-stresses, as well as a $B_1$- and a $B_2$-symmetric self-stress (that is, an anti-symmetric self-stress for each pair of perpendicular mirrors) in this case, which is a useful property for the construction of gridshells. 

Here we chose $k=-9$ and  $e_{\sigma_v}=e_{\sigma_d}=8$ to meet these conditions. See the  $\mathcal{C}_{4v}$ count below for full details on the detected self-stresses and mechanisms. The counts below also show that using increasingly large symmetry groups strictly increases the number of symmetry-detectable self-stresses (and mechanisms) for this example:
\begin{itemize}
\item $\mathcal{C}_s$: $\Gamma(m)-\Gamma(s)= -8A'-A''$, so we find 9 self-stresses;
\item $\mathcal{C}_{2v}$: $\Gamma(m)-\Gamma(s)= -6A_1+A_2-2B_1-2B_2$, so we find 10 self-stresses;
\item $\mathcal{C}_{4v}$: $\Gamma(m)-\Gamma(s)= -5A_1+2A_2-2B_1-2B_2-2E$, so we find 11 self-stresses.
\end{itemize}
Note that since infinitesimal rigidity is projectively invariant, we may use projective transformations to reduce the  $\mathcal{C}_{4v}$ symmetry to a desired  subgroup while preserving the dimension of the space of self-stresses. The framework in Figure~\ref{fig:spiderweb}(b), for example, is obtained from the one in (a) by an affine transformation, and so we know from the $\mathcal{C}_{4v}$ analysis that it must also have at least 11 self-stresses.  Such an analysis of a projectively equivalent framework with a larger symmetry group can be a useful tool for finding additional self-stresses.
\end{example}

\section{Methods beyond symmetry} \label{sec:nonsym}

\begin{figure}[htp]
\begin{center}
\begin{tikzpicture}[very thick,scale=0.7]
\tikzstyle{every node}=[circle, draw=black, fill=white, inner sep=0pt, minimum width=3pt];
   
       \path (90:2.2cm) node (p1)  {} ;
     \path (162:2.2cm) node (p2)  {} ;
     
      \path (18:2.2cm) node (p5)  {} ;
      
     \path (0,-0.15) node (p6)  {} ;
     \path (0,0.42) node (p7)  {} ;
      \path (-1,0.7) node (p8)  {} ;
     \path (1,0.7) node (p9)  {} ;

       \node [draw=white, fill=white] (b) at (0.25,2.5) {$b$};
       \node [draw=white, fill=white] (b) at (-6.3,2.5) {$a$};
       \node [draw=white, fill=white] (b) at (6.3,2.5) {$c$};

     \draw[dashed,red,thin] (-6.5,2.2)--(6.5,2.2);
      \draw[dashed,red,thin] (p7)--(6,2.2);
      \draw[dashed,red,thin] (p6)--(6,2.2);
        \draw[dashed,red,thin] (p7)--(-6,2.2);
      \draw[dashed,red,thin] (p6)--(-6,2.2);

     \draw(p1)--(p2);
   \draw(p1)--(p8);
   \draw(p1)--(p9);
   \draw(p1)--(p5);
     \draw(p8)--(p2);
   \draw(p5)--(p9);
   \draw(p6)--(p7);
   \draw(p6)--(p2);
    \draw(p6)--(p5);
       \draw(p8)--(p7);
    \draw(p7)--(p9);

    \draw[dashed, thin](0,-0.8)--(0,2.7);

        \end{tikzpicture}
      \end{center}
\vspace{-0.6cm} \caption{A framework with the same underlying graph as in Figure~\ref{fig:symfw}(a) satisfying the pure condition for this graph: the points $a, b$ and $c$ are collinear.} \label{fig:pure}
\end{figure}
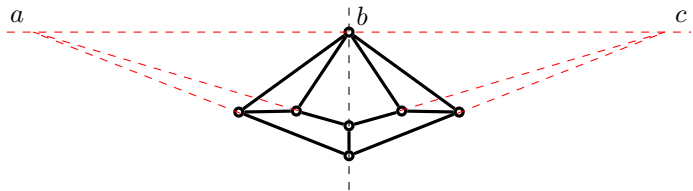

While the symmetry-based method presented in this paper  provides a useful tool for increasing the number of independent states of self-stress in frameworks, it does not,  in general, find the maximum possible number of independent self-stresses for a given graph and  symmetry group.  This is because the existence of self-stresses  is a \emph{projective} geometric condition, and not a symmetric condition. Consider, for example, the framework in Figure~\ref{fig:symfw}(a). This framework has $k=0$ and a symmetry analysis with the point group $\mathcal{C}_s$ detects no self-stress or mechanism. In fact, since $e_\sigma=1$,  we obtain an isostatic framework for all `generic' positions of the vertices (i.e., almost all positions of the vertices satisfying the reflection symmetry constraint), as shown in \cite{BS4}. However, if the vertices are placed in a special geometric position satisfying the so-called \emph{pure condition} of the graph (see \cite[Table 1]{WW} and Figure~\ref{fig:pure}), then the framework has a non-trivial self-stress and mechanism.  Note that both the self-stress and the mechanism are fully-symmetric so that $\Gamma(m)-\Gamma(s)=(0,0)=A'-A'$,  and hence they are not detected with the symmetry-extended Maxwell rule.

The pure conditions for some small standard graphs are well known (see \cite[Table 1]{WW}, for example). In general, however, finding the special geometric conditions which give rise to additional self-stresses that are not detected with the symmetry-extended Maxwell count requires a  non-trivial analysis.

\begin{figure}[htp]
\begin{center}
\begin{tikzpicture}[very thick,scale=0.6]
\tikzstyle{every node}=[circle, draw=black, fill=white, inner sep=0pt, minimum width=4pt];
   
       \path (0,0) node (p1)  {} ;
    \path (0,3) node (p2)  {} ;
    \path (-0.8,2.3) node (p3)  {} ;
     \path (-1.3,0.7) node (p4)  {} ;
      \path (-1.8,2.5) node (p5)  {} ;
     \path (-1.8,0.2) node (p6)  {} ;
    
    \path (0.8,2.3) node (p3r)  {} ;
     \path (1.3,0.7) node (p4r)  {} ;
      \path (1.8,2.5) node (p5r)  {} ;
     \path (1.8,0.2) node (p6r)  {} ;
     
               \draw (p1)  --  (p2);
      \draw (p1)  --  (p4);
        \draw (p1)  --  (p6);
     \draw (p2)  --  (p5);
        \draw (p2)  --  (p3);
      \draw (p3)  --  (p5);
        \draw (p4)  --  (p6);
               \draw (p3)  --  (p4);
              \draw (p5)  --  (p6);
              
                \draw (p1)  --  (p4r);
        \draw (p1)  --  (p6r);
     \draw (p2)  --  (p5r);
        \draw (p2)  --  (p3r);
      \draw (p3r)  --  (p5r);
        \draw (p4r)  --  (p6r);
               \draw (p3r)  --  (p4r);
              \draw (p5r)  --  (p6r);
                         
    \draw[dashed, thin](0,-0.5)--(0,3.6); 
      
 \node [draw=white, fill=white] (b) at (0,-0.9) {(a)};
        \end{tikzpicture}
        \hspace{1.5cm}
      \begin{tikzpicture}[very thick,scale=0.6]
\tikzstyle{every node}=[circle, draw=black, fill=white, inner sep=0pt, minimum width=4pt];
      \path (0,0) node (p1)  {} ;
    \path (0,3) node (p2)  {} ;
    \path (-1,2.3) node (p3)  {} ;
     \path (-1,0.7) node (p4)  {} ;
      \path (-1.8,2.5) node (p5)  {} ;
     \path (-1.8,0.2) node (p6)  {} ;
    
    \path (1,2.3) node (p3r)  {} ;
     \path (1,0.7) node (p4r)  {} ;
      \path (1.8,2.5) node (p5r)  {} ;
     \path (1.8,0.2) node (p6r)  {} ;
     
               \draw (p1)  --  (p2);
      \draw (p1)  --  (p4);
        \draw (p1)  --  (p6);
     \draw (p2)  --  (p5);
        \draw (p2)  --  (p3);
      \draw (p3)  --  (p5);
        \draw (p4)  --  (p6);
               \draw (p3)  --  (p4);
              \draw (p5)  --  (p6);
              
                \draw (p1)  --  (p4r);
        \draw (p1)  --  (p6r);
     \draw (p2)  --  (p5r);
        \draw (p2)  --  (p3r);
      \draw (p3r)  --  (p5r);
        \draw (p4r)  --  (p6r);
               \draw (p3r)  --  (p4r);
              \draw (p5r)  --  (p6r);
                         
    \draw[dashed, thin](0,-0.5)--(0,3.6); 
      
 \node [draw=white, fill=white] (b) at (0,-0.9) {(b)};
        \end{tikzpicture}
\end{center}
\vspace{-0.6cm} \caption{Two frameworks with $\mathcal{C}_s$ symmetry.  (a) is isostatic, but  (b) has a fully-symmetric and an anti-symmetric self-stress (and two corresponding mechanisms) since the triangular prism subgraph on either side of the mirror is placed in a special position satisfying the pure condition for this graph: each of the two frameworks forms a Desargues configuration \cite[Table 1]{WW}.} \label{fig:des}
\end{figure}
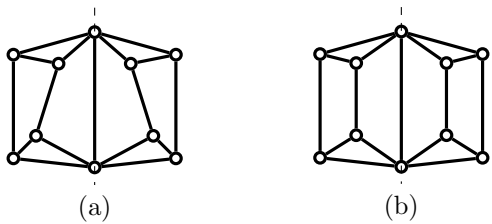

Given a framework with a reflection symmetry, it is natural to try to create further self-stresses -- in addition to the ones detected by the symmetry-extended Maxwell rule -- by placing the vertices on one side of the mirror in a special position so that this part of the structure becomes self-stressed. This self-stress is then duplicated on the other side of the mirror, creating a fully-symmetric and an anti-symmetric self-stress for the whole framework. None of these self-stresses can be detected with the method presented in this paper since they are created independently from the reflection symmetry.  See Figure~\ref{fig:des} for an example.

The example in Figure~\ref{fig:des} suggests that we may `glue together' self-stressed frameworks to build up larger frameworks with many independent states of self-stress. Note, however, that this method of gluing together  framework primitives is problematic from a practical point of view. One of Maxwell's seminal papers \cite{Maxwell1864} states that if a planar framework  possesses a state of self-stress then it must be the vertical projection of a plane-faced polyhedron (which is also known as the \emph{discrete Airy stress function polyhedron}). Sometimes a vertical lifting of the form diagram is taken as a gridshell roof, since this guarantees planarity of faces which has beneficial properties in terms of cost and construction. By gluing together framework primitives, the edge of each primitive often remains on the $z=0$ plane for each lifting and this is not architecturally acceptable. 

Similarly, we may start with a planar framework and subdivide its faces -- either by inserting additional bars or by inserting entire self-stressed frameworks -- to create further states of self-stress. If we simply insert additional bars, then this of course also decreases the freedom number $k$ of the framework. It is possible to insert self-stressed frameworks into the faces of a given planar framework without changing its freedom number,  but this method of subdividing faces has the same practical problems as gluing framework primitives together, since the newly created self-stresses are all local.

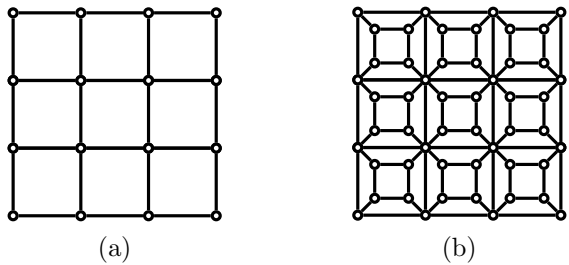
\begin{figure}[htp]
\begin{center}
\begin{tikzpicture}[very thick,scale=0.45]
\tikzstyle{every node}=[circle, draw=black, fill=white, inner sep=0pt, minimum width=3pt];
    
       \path (-1,0) node (p11)  {} ;
    \path (1,0) node (p22)  {} ;
    \path (-1,-2) node (p33)  {} ;
     \path (1,-2) node (p44)  {} ;
     
      \path (3,0) node (p1r)  {} ;
    \path (3,-2) node (p2r)  {} ;
    \path (5,0) node (p3r)  {} ;
     \path (5,-2) node (p4r)  {} ;

        \draw (p11)  --  (p22);
      \draw (p11)  --  (p33);
        \draw (p22)  --  (p44);
               \draw (p33)  --  (p44);

 \draw (p22)  --  (p1r);
               \draw (p44)  --  (p2r);   
                 \draw (p1r)  --  (p3r);
                  \draw (p2r)  --  (p4r);
  \draw (p2r)  --  (p1r);
               \draw (p4r)  --  (p3r);   
               
            \path (-1,-2) node (b11)  {} ;
    \path (1,-2) node (b22)  {} ;
    \path (-1,-4) node (b33)  {} ;
     \path (1,-4) node (b44)  {} ;
     
      \path (3,-2) node (b1r)  {} ;
    \path (3,-4) node (b2r)  {} ;
    \path (5,-2) node (b3r)  {} ;
     \path (5,-4) node (b4r)  {} ;

        \draw (b11)  --  (b22);
      \draw (b11)  --  (b33);
        \draw (b22)  --  (b44);
               \draw (b33)  --  (b44);

 \draw (b22)  --  (b1r);
               \draw (b44)  --  (b2r);   
                 \draw (b1r)  --  (b3r);
                  \draw (b2r)  --  (b4r);
  \draw (b2r)  --  (b1r);
               \draw (b4r)  --  (b3r);     
              
            \path (-1,-4) node (c11)  {} ;
    \path (1,-4) node (c22)  {} ;
    \path (-1,-6) node (c33)  {} ;
     \path (1,-6) node (c44)  {} ;
     
      \path (3,-4) node (c1r)  {} ;
    \path (3,-6) node (c2r)  {} ;
    \path (5,-4) node (c3r)  {} ;
     \path (5,-6) node (c4r)  {} ;

        \draw (c11)  --  (c22);
      \draw (c11)  --  (c33);
        \draw (c22)  --  (c44);
               \draw (c33)  --  (c44);

 \draw (c22)  --  (c1r);
               \draw (c44)  --  (c2r);   
                 \draw (c1r)  --  (c3r);
                  \draw (c2r)  --  (c4r);
  \draw (c2r)  --  (c1r);
               \draw (c4r)  --  (c3r);

 \node [draw=white, fill=white] (b) at (2,-7) {(a)};
        \end{tikzpicture}
        \hspace{1.5cm}
      \begin{tikzpicture}[very thick,scale=0.45]
\tikzstyle{every node}=[circle, draw=black, fill=white, inner sep=0pt, minimum width=3pt];
   
  \path (-0.5,-0.5) node (p1)  {} ;
    \path (0.5,-0.5) node (p2)  {} ;
    \path (-0.5,-1.5) node (p3)  {} ;
     \path (0.5,-1.5) node (p4)  {} ;
     
        \path (-1,0) node (p11)  {} ;
    \path (1,0) node (p22)  {} ;
    \path (-1,-2) node (p33)  {} ;
     \path (1,-2) node (p44)  {} ;
     
      \path (3,0) node (p1r)  {} ;
    \path (3,-2) node (p2r)  {} ;
    \path (5,0) node (p3r)  {} ;
     \path (5,-2) node (p4r)  {} ;

      \path (1.5,-0.5) node (a1)  {} ;
    \path (2.5,-0.5) node (a2)  {} ;
    \path (1.5,-1.5) node (a3)  {} ;
     \path (2.5,-1.5) node (a4)  {} ;
     
      \path (3.5,-0.5) node (b1)  {} ;
    \path (4.5,-0.5) node (b2)  {} ;
    \path (3.5,-1.5) node (b3)  {} ;
     \path (4.5,-1.5) node (b4)  {} ;
     
          \draw (a1)  --  (a2);
      \draw (a1)  --  (a3);
        \draw (a2)  --  (a4);
     \draw (a3)  --  (a4);
      
        \draw (a1)  --  (p22);
      \draw (a3)  --  (p44);
        \draw (a2)  --  (p1r);
               \draw (a4)  --  (p2r);
     
       \draw (b1)  --  (b2);
      \draw (b1)  --  (b3);
        \draw (b2)  --  (b4);
     \draw (b3)  --  (b4);
      
        \draw (b1)  --  (p1r);
      \draw (b3)  --  (p2r);
        \draw (b2)  --  (p3r);
               \draw (b4)  --  (p4r);

               \draw (p1)  --  (p2);
      \draw (p1)  --  (p3);
        \draw (p2)  --  (p4);
     \draw (p3)  --  (p4);
      
        \draw (p11)  --  (p22);
      \draw (p11)  --  (p33);
        \draw (p22)  --  (p44);
               \draw (p33)  --  (p44);
            
              \draw (p1)  --  (p11);
               \draw (p2)  --  (p22);   
                 \draw (p3)  --  (p33);
                  \draw (p4)  --  (p44);
                  
                   \draw (p22)  --  (p1r);
               \draw (p44)  --  (p2r);   
                 \draw (p1r)  --  (p3r);
                  \draw (p2r)  --  (p4r);
  \draw (p2r)  --  (p1r);
               \draw (p4r)  --  (p3r); 
      
      
       \path (-0.5,-2.5) node (p1y)  {} ;
    \path (0.5,-2.5) node (p2y)  {} ;
    \path (-0.5,-3.5) node (p3y)  {} ;
     \path (0.5,-3.5) node (p4y)  {} ;
     
        \path (-1,-2) node (p11y)  {} ;
    \path (1,-2) node (p22y)  {} ;
    \path (-1,-4) node (p33y)  {} ;
     \path (1,-4) node (p44y)  {} ;
     
      \path (3,-2) node (p1ry)  {} ;
    \path (3,-4) node (p2ry)  {} ;
    \path (5,-2) node (p3ry)  {} ;
     \path (5,-4) node (p4ry)  {} ;

      \path (1.5,-2.5) node (a1y)  {} ;
    \path (2.5,-2.5) node (a2y)  {} ;
    \path (1.5,-3.5) node (a3y)  {} ;
     \path (2.5,-3.5) node (a4y)  {} ;
     
      \path (3.5,-2.5) node (b1y)  {} ;
    \path (4.5,-2.5) node (b2y)  {} ;
    \path (3.5,-3.5) node (b3y)  {} ;
     \path (4.5,-3.5) node (b4y)  {} ;
     
          \draw (a1y)  --  (a2y);
      \draw (a1y)  --  (a3y);
        \draw (a2y)  --  (a4y);
     \draw (a3y)  --  (a4y);
      
        \draw (a1y)  --  (p22y);
      \draw (a3y)  --  (p44y);
        \draw (a2y)  --  (p1ry);
               \draw (a4y)  --  (p2ry);
     
       \draw (b1y)  --  (b2y);
      \draw (b1y)  --  (b3y);
        \draw (b2y)  --  (b4y);
     \draw (b3y)  --  (b4y);
      
        \draw (b1y)  --  (p1ry);
      \draw (b3y)  --  (p2ry);
        \draw (b2y)  --  (p3ry);
               \draw (b4y)  --  (p4ry);

               \draw (p1y)  --  (p2y);
      \draw (p1y)  --  (p3y);
        \draw (p2y)  --  (p4y);
     \draw (p3y)  --  (p4y);
      
        \draw (p11y)  --  (p22y);
      \draw (p11y)  --  (p33y);
        \draw (p22y)  --  (p44y);
               \draw (p33y)  --  (p44y);
            
              \draw (p1y)  --  (p11y);
               \draw (p2y)  --  (p22y);   
                 \draw (p3y)  --  (p33y);
                  \draw (p4y)  --  (p44y);
                  
                   \draw (p22y)  --  (p1ry);
               \draw (p44y)  --  (p2ry);   
                 \draw (p1ry)  --  (p3ry);
                  \draw (p2ry)  --  (p4ry);
  \draw (p2ry)  --  (p1ry);
               \draw (p4ry)  --  (p3ry); 
      
       \path (-0.5,-4.5) node (p1w)  {} ;
    \path (0.5,-4.5) node (p2w)  {} ;
    \path (-0.5,-5.5) node (p3w)  {} ;
     \path (0.5,-5.5) node (p4w)  {} ;
     
        \path (-1,-4) node (p11w)  {} ;
    \path (1,-4) node (p22w)  {} ;
    \path (-1,-6) node (p33w)  {} ;
     \path (1,-6) node (p44w)  {} ;
     
      \path (3,-4) node (p1rw)  {} ;
    \path (3,-6) node (p2rw)  {} ;
    \path (5,-4) node (p3rw)  {} ;
     \path (5,-6) node (p4rw)  {} ;

      \path (1.5,-4.5) node (a1w)  {} ;
    \path (2.5,-4.5) node (a2w)  {} ;
    \path (1.5,-5.5) node (a3w)  {} ;
     \path (2.5,-5.5) node (a4w)  {} ;
     
      \path (3.5,-4.5) node (b1w)  {} ;
    \path (4.5,-4.5) node (b2w)  {} ;
    \path (3.5,-5.5) node (b3w)  {} ;
     \path (4.5,-5.5) node (b4w)  {} ;
     
          \draw (a1w)  --  (a2w);
      \draw (a1w)  --  (a3w);
        \draw (a2w)  --  (a4w);
     \draw (a3w)  --  (a4w);
      
        \draw (a1w)  --  (p22w);
      \draw (a3w)  --  (p44w);
        \draw (a2w)  --  (p1rw);
               \draw (a4w)  --  (p2rw);
     
       \draw (b1w)  --  (b2w);
      \draw (b1w)  --  (b3w);
        \draw (b2w)  --  (b4w);
     \draw (b3w)  --  (b4w);
      
        \draw (b1w)  --  (p1rw);
      \draw (b3w)  --  (p2rw);
        \draw (b2w)  --  (p3rw);
               \draw (b4w)  --  (p4rw);

               \draw (p1w)  --  (p2w);
      \draw (p1w)  --  (p3w);
        \draw (p2w)  --  (p4w);
     \draw (p3w)  --  (p4w);
      
        \draw (p11w)  --  (p22w);
      \draw (p11w)  --  (p33w);
        \draw (p22w)  --  (p44w);
               \draw (p33w)  --  (p44w);
            
              \draw (p1w)  --  (p11w);
               \draw (p2w)  --  (p22w);   
                 \draw (p3w)  --  (p33w);
                  \draw (p4w)  --  (p44w);
                  
                   \draw (p22w)  --  (p1rw);
               \draw (p44w)  --  (p2rw);   
                 \draw (p1rw)  --  (p3rw);
                  \draw (p2rw)  --  (p4rw);
  \draw (p2rw)  --  (p1rw);
               \draw (p4rw)  --  (p3rw);
 \node [draw=white, fill=white] (b) at (2,-7) {(b)};
        \end{tikzpicture}
\end{center}
\vspace{-0.6cm} \caption{(a) A framework with freedom number $k=5$ which has no self-stress. The framework in (b) also has $k=5$ and is obtained by inserting a self-stressed framework into each face of (a). So it has 9 independent self-stresses, one for each original face.} \label{fig:quadrefined}
\end{figure}

Consider, for example, the planar framework shown in Figure~\ref{fig:quadrefined}(a). If we subdivide each of the quadrilateral faces by inserting a cube graph (see Figure~\ref{fig:symfw}(b)), then the framework is kept quad-dominant and the freedom number remains unchanged. Moreover, by placing the newly added vertices in suitable geometric positions, the aspect ratio of the quadrilaterals is kept within an acceptable range, and an independent self-stress is created within each original face (see Figure~\ref{fig:quadrefined}(b) and recall Figure~\ref{fig:symfw}(b)). However, since the self-stresses in this refined framework are all local, its vertical lifting may not yield a suitable structure for a gridshell roof.

\section{Pinned frameworks} \label{sec:pinned}

All of the above immediately transfers to pinned frameworks, where the rigid body motions have been eliminated by the pinning of some vertices. For a \emph{pinned} bar-joint framework in the plane, the Maxwell rule becomes
\begin{equation}
    m-s=2v-e,
\label{eq:scrulepinned}
\end{equation}
where $v$ is the number of \emph{internal} (or \emph{unpinned}) vertices. 
Similarly, as shown in \cite{FGsymmax}, the symmetry-extended counting rule for pinned  frameworks simplifies to
\begin{equation}
    \Gamma(m)-\Gamma(s)=\Gamma(v)\times \Gamma_{\mathrm{T}}-\Gamma(e).
\label{eq:bj1pinned}\end{equation}

For pinned frameworks, the character calculations given in Table~\ref{tab:2D} simplify as shown in Table~\ref{tab:2Dpinned}.
\begin{table}[ht]\begin{center}
    \begin{tabular}{l|c c c c}
                        & $\phantom{-}E$       &  $C_{n > 2}$   &
                            $\phantom{-}C_2$   & $\phantom{-}\sigma$  \\ \hline
    $\phantom{=}\Gamma(v)$& $\phantom{-}v$     & $v_c$             &
                            $\phantom{-}v_c$   & $\phantom{-}v_\sigma$\\
    $\phantom{=}\times \Gamma_{T}$&   $\phantom{-}2$     & $2\cos\phi$       &
                            $-2$    & $\phantom{-}0$       \\ \hline
    $=\Gamma(v)\times\Gamma_{T}$
                        & $\phantom{-}2v$      & $2 v_c \cos\phi$  &
                            $-2 v_c$& $\phantom{-}0$       \\
    $\phantom{=}-\Gamma(e)$        & $-e$      & $0$               &
                            $-e_2$  & $-e_\sigma$\\\hline
       $=\Gamma(m) - \Gamma(s)$
                        & $2v-e$  & $2v_c\cos\phi$ &
                            $-2v_c - e_2$ & $-e_\sigma$

    \end{tabular}
    \caption{Calculations of characters for the 2D
        symmetry-extended Maxwell equation for pinned frameworks (\ref{eq:bj1pinned}). Note the similarity to Table~\ref{tab:2D}.}
    \label{tab:2Dpinned}
    \end{center}
\end{table}

In the following we will consider planar pinned frameworks satisfying the count $m-s=2v-e=k$. As before, we will call the integer $k$ the \emph{freedom number} of the framework.
We may obtain formulas for creating states of self-stress (or mechanisms) in symmetric pinned frameworks in the analogous way as for unpinned frameworks. We summarise the formulas for some basic groups  below.

\begin{itemize}
\item For a framework with reflection  symmetry $\mathcal{C}_{s}$, we obtain:
$$\Gamma(m)-\Gamma(s)= (k,-e_\sigma)=\frac{k-e_\sigma}{2}A' + \frac{k+e_\sigma}{2}A''.$$

\item For a framework with half-turn  symmetry $\mathcal{C}_{2}$, we obtain:
$$\Gamma(m)-\Gamma(s)= (k,-2v_c-e_2)=\begin{cases}
\frac{k}{2}A + \frac{k}{2}B, &  \textrm{ if } v_c=e_2=0\\
\frac{k-1}{2}A + \frac{k+1}{2}B, &  \textrm{ if } v_c=0, e_2=1\\
\frac{k-2}{2}A + \frac{k+2}{2}B, &  \textrm{ if } v_c=1\\
\end{cases} $$

\item For a framework with dihedral symmetry $\mathcal{C}_{2v}$, we obtain:
$$\Gamma(m)-\Gamma(s)= (k,-2v_c-e_2,-e_{\sigma_h},-e_{\sigma_v}),$$
which leads to the following formulas for $\Gamma(m)-\Gamma(s)$.
\end{itemize}

For $v_c=0$ and $e_2=0$ we obtain: 
{
\medmuskip=0mu
\thinmuskip=0mu
\thickmuskip=0mu
\nulldelimiterspace=0pt
\scriptspace=0pt
\begin{equation*}
\Gamma(m)-\Gamma(s)=\frac{k-e_{\sigma_h}-e_{\sigma_v}}{4}A_1 + \frac{k+e_{\sigma_h}+e_{\sigma_v}}{4}A_2 +\frac{k-e_{\sigma_h}+e_{\sigma_v}}{4}B_1 + \frac{k+e_{\sigma_h}-e_{\sigma_v}}{4}B_2 \label{eq:pinned_v0}\end{equation*}
}

For $v_c=0$ and $e_2=1$ we obtain: 
{
\medmuskip=0mu
\thinmuskip=0mu
\thickmuskip=0mu
\nulldelimiterspace=0pt
\scriptspace=0pt
\begin{equation*}
\Gamma(m)-\Gamma(s)=\frac{k-e_{\sigma_h}-e_{\sigma_v}-1}{4}A_1 + \frac{k+e_{\sigma_h}+e_{\sigma_v}-1}{4}A_2 +\frac{k-e_{\sigma_h}+e_{\sigma_v}+1}{4}B_1 + \frac{k+e_{\sigma_h}-e_{\sigma_v}+1}{4}B_2 \label{eq:pinned_v0}\end{equation*}
}

For  $v_c=1$ we obtain:
{
\medmuskip=0mu
\thinmuskip=0mu
\thickmuskip=0mu
\nulldelimiterspace=0pt
\scriptspace=0pt
\begin{equation*}\Gamma(m)-\Gamma(s)=\frac{k-e_{\sigma_h}-e_{\sigma_v}-2}{4}A_1 + \frac{k+e_{\sigma_h}+e_{\sigma_v}-2}{4}A_2 +\frac{k-e_{\sigma_h}+e_{\sigma_v}+2}{4}B_1 + \frac{k+ e_{\sigma_h}-e_{\sigma_v}+2}{4}B_2\end{equation*}
}

In each case, we can draw analogous conclusions regarding the states of self-stress of the framework  as for unpinned frameworks above. We leave this discussion, as well as the straightforward computations for other groups, to the reader. We conclude this section with a practical example instead.

\begin{figure}[htp]\begin{center} 
      \includegraphics[scale=0.68]{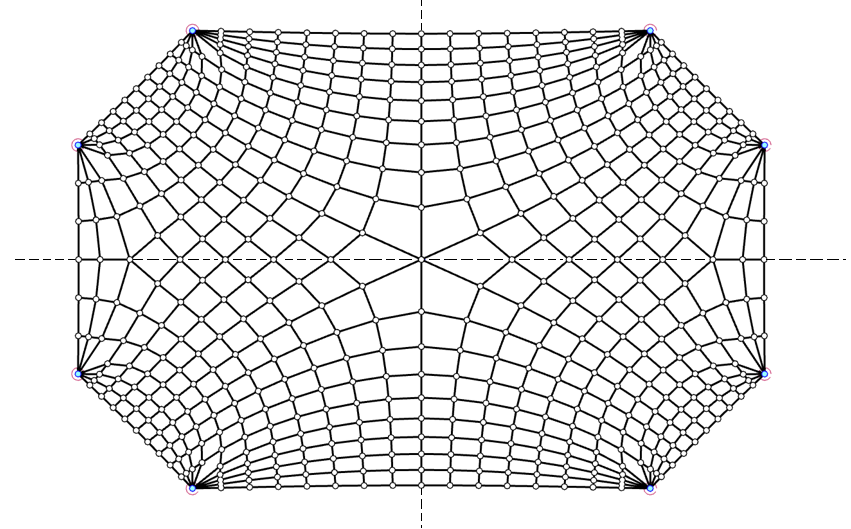}\end{center}
      \caption{The form diagram of the gridshell structure in Figure~\ref{fig:pic}. As a pinned framework it has  freedom number $k=4$ and $\mathcal{C}_{2v}$ symmetry. (The pinned vertices are shown in blue.) This framework has 7 independent symmetry-detectable self-stresses, 5 of which are fully-symmetric and 2 of which are anti-symmetric with respect to $\sigma_h$ (the reflection in the horizontal mirror). }
\label{fig:grid1}
\end{figure}

\begin{example} The pinned framework in Figure~\ref{fig:grid1} is the form diagram of the gridshell structure shown in Figure~\ref{fig:pic} (i.e., it is  the vertical projection of the gridhell structure onto the $xy$-plane).  This framework has $e=2v-4=1102$, 
so $k=4$. It was form-found using the force density method \cite{Schek} so it is known to have at least one fully-symmetric self-stress. A symmetry analysis reveals significant further information. The framework has $v_c=1$, $e_{\sigma_v}=18$  and $e_{\sigma_h}=4$. So if we analyse it with the full $\mathcal{C}_{2v}$  symmetry, then we obtain 
$$\Gamma(m)-\Gamma(s)= -5A_1 + 6A_2 +5B_1 - 2B_2.$$
Thus, we see that this framework has at least 5 self-stresses that are fully-symmetric and 2 self-stresses of symmetry $B_2$. (We also detect 6 mechanisms of symmetry $A_2$ and 5 mechanisms of symmetry $B_1$.) As has previously been discussed (see Section~\ref{sec:c2v} and note that the reasoning is analogous for unpinned and pinned frameworks), the existence of the $B_2$-symmetric self-stresses is a consequence of the large difference between $e_{\sigma_v}$ and $e_{\sigma_h}$.

Note that an analysis of  the framework with $\mathcal{C}_s$ symmetry, where the reflection is in the vertical mirror, also finds $7$  self-stresses all of which are fully-symmetric with respect to the vertical mirror.  A $\mathcal{C}_s$ analysis with the other reflection does not find any self-stresses.
\end{example}

\section{Further comments and future work} \label{sec:future}

 The methods of this paper can easily be extended to non-planar frameworks. However, since for non-planar frameworks there can be multiple bars that are unshifted by a $C_2$ rotation, for example, some of the formulas become slightly more involved. The methods can also be extended to frameworks in 3-space, which has potential applications in the analysis of space frames, for example. 
 
As has previously been discussed, this paper provides an efficient method for increasing the number of independent states of self-stress in symmetric frameworks. However, finding a realisation of a given graph that has the maximum possible number of states of self-stress (with or without a specified point group symmetry) remains a challenging open problem. Note that maximising the space of self-stresses is equivalent to maximising the space of mechanisms or parallel redrawings (see \cite{SW1,Wmatroids}, for example), or the decomposibility of the discrete Airy stress function polyhedron \cite{Smilansky}, so there are several different but equivalent ways to formulate this problem. 

An important tool in analysing a form diagram is the \emph{reciprocal diagram} or \emph{force diagram}, which is a geometric construction that has appeared, independently, in areas such as graphical statics, rigidity theory, scene analysis and computational geometry since the time of Maxwell \cite{SW1}. In a recent paper McRobie et al describe the relationship between mechanisms and states of self-stress in the form and force diagrams \cite{Allan2015}. 
 It would be interesting to investigate this relationship with an emphasis on symmetry. This is left to a future paper.

Finally, it would be useful to establish procedures for subdividing faces of a planar framework in such a way that  additional non-local self-stresses are created. This is left as another area of future research.

\section*{Acknowledgements}

We thank the Fields Institute for hosting the 2021 workshop on ‘Progress and Open Problems in Rigidity Theory’ during which this work was started. We also thank Allan McRobie and Toby Mitchell for helpful discussions.

\section*{Appendix} \label{app}

The expression~(\ref{eq:cn11}) for $\Gamma(m)-\Gamma(s)$ for frameworks with rotational symmetry follows from the following proposition.

\begin{proposition} \label{prop:triv} For $t\in\{1,\ldots, n-1\}$ and $\epsilon=e^{\frac{2\pi i}{n}}$, we have $$\sum_{j=0}^{n-1} \epsilon^{tj}\cos \Big(\frac{j2\pi}{n}\Big)=\begin{cases}
\frac{n}{2} &  \textrm{ if } t=1 \textrm{ or } n-1 \\
0 &  \textrm{ otherwise}\\
\end{cases}$$
\end{proposition}
\emph{Proof.} It is well known that the sum of the entries of each character $A_t$, $\sum_{j=0}^{n-1}\epsilon^{tj}$, is zero for each $t\in\{1,\ldots, n-1\}$. Suppose first that $t\in\{2,\ldots, n-2\}$. From the trigonometric identities $\cos x \cos y = \frac{1}{2}\big(\cos(x-y)+\cos (x+y)\big)$ and $\sin x \cos y = \frac{1}{2}\big(\sin(x+y)+\sin (x-y)\big)$ we obtain for $\sum_{j=0}^{n-1} \epsilon^{tj}\cos \Big(\frac{j2\pi}{n}\Big)$:
\begin{eqnarray*} && \sum_{j=0}^{n-1}\Big(\cos\Big(\frac{tj2\pi}{n}\Big) +i\sin\Big(\frac{tj2\pi}{n}\Big) \Big)\cos\Big(\frac{j2\pi}{n}\Big)\\
&=&  \frac{1}{2}\sum_{j=0}^{n-1}\cos\Big(\frac{(t-1)j2\pi }{n}\Big) +\cos\Big(\frac{(t+1)j2\pi}{n}\Big)+ i\Big(\sin\Big(\frac{(t+1)j2\pi}{n}\Big) +\sin\Big(\frac{(t-1)j2\pi }{n}  \Big) \Big)\\
&=& \frac{1}{2}\sum_{j=0}^{n-1}\Big(\epsilon^{(t-1)j} + \epsilon^{(t+1)j}\Big)\\
& = & 0\end{eqnarray*}
since $1\leq t-1<t+1\leq n-1$.

Suppose next that $t=1$. Then $$\sum_{j=0}^{n-1} \epsilon^{j}\cos \Big(\frac{j2\pi}{n}\Big)=  \sum_{j=0}^{n-1}\Big(\cos^2\Big(\frac{j2\pi}{n}\Big) +i\sin\Big(\frac{j2\pi}{n}\Big) \cos\Big(\frac{j2\pi}{n} \Big)\Big).$$
Now, using the trigonometric identity $\cos^2 x=\frac{1}{2}\cos 2x+1$, we have 
$$\sum_{j=0}^{n-1}\cos^2\Big(\frac{j2\pi}{n}\Big)=\frac{1}{2}\sum_{j=0}^{n-1}\Big(\cos\Big(\frac{j4\pi}{n}\Big)+1\Big) =\frac{n}{2}$$
since $$\sum_{j=0}^{n-1}\cos\Big(\frac{j4\pi}{n}\Big)= \textrm{Re} \sum_{j=0}^{n-1} \epsilon^{2j}=\textrm{Re}\Big(\frac{1- \epsilon^{2n}}{1- \epsilon^{2}}\Big)=0.$$
Also, using the  trigonometric identity $\sin x\cos x=\frac{1}{2}\sin 2x$, we have 
$$\sum_{j=0}^{n-1}i\sin\Big(\frac{j2\pi}{n}\Big) \cos\Big(\frac{j2\pi}{n} \Big)=\frac{i}{2}\sum_{j=0}^{n-1}\sin\Big(\frac{j4\pi}{n}\Big)=\frac{i}{2}\textrm{Im}\Big(\sum_{j=0}^{n-1} \epsilon^{2j}\Big)=0.$$
Finally, if $t=n-1$, then the result follows from the argument for $t=1$ and the fact that cosine and sine are even and odd functions, respectively.  \hfill $\square$

\end{document}